\documentstyle[12pt]{article}

\begin{document}
\title{Wrap Groups of Non-Archimedean Fiber Bundles.}

\author{S. V. Ludkovsky.}
\date{15 September 2011 \thanks{Mathematics subject classification
 22E20, 46S10, 54H15 and 57S20.} \thanks{keywords:
wrap group; fiber bundle; infinite field; non-archimedean norm;
Cayley-Dickson algebra}} \maketitle

\begin{abstract}
\par Fiber bundles over infinite fields with non-trivial
ultra-norms are considered. For them geometric wrap groups are
defined and investigated. Besides fields also Cayley-Dickson
algebras over fields of characteristic not equal to two are taken
into account. For fibers over them wrap groups are introduced and
their structure is investigated. Different classes of smoothness for
wrap groups are used. It is demonstrated that generally such groups
are infinite dimensional over the corresponding field and totally
disconnected groups. That is, they are continuous or differentiable
non-archimedean differentiable uniform spaces and the composition
$(f,g)\mapsto f^{-1}g$ is continuous or differentiable depending on
a class of smoothness of groups. Skew products of wrap groups are
studied as well.
\end{abstract}

\section{Introduction.}
\par Groups of geometric loops and paths for real manifolds were introduced
in the 1930-th and they are very important in differential geometry,
algebraic topology and theoretical physics \cite{bryl,mensk,swit}.
Possibly the first author was S. Lefshetz who studied them over real
manifolds. He used families of continuous mappings, that was rather
restrictive and led him to the necessity to combine a geometric
construction with an additional algebraic construction with some
elements of free groups. Then more natural approach was proposed by
J. Milnor in the class of Sobolev mappings. Later on they were
generalized for real and complex fiber bundles \cite{bryl}.
Previously loop groups were considered as classes of mappings from
the unit circle $S^1$ into a real or complex fiber bundle with a
parallel transport structure. For spheres iterated loop groups were
considered using reduced products of copies of circles. Recently
they were generalized as wrap groups for rather common manifolds and
fibers over the real and complex fields, the quaternion skew field
and the octonion algebra \cite{lugmlg,luijmgta09,ludan,lujmslg},
where also some new structural theorems were proved. In the latter
work Sobolev classes of smoothness were used. For general manifolds
different from circles and spheres earlier geometric interpretation
is already lost, so such groups were called wrap groups.
\par For manifolds over non-archimedean fields of zero characteristic
loop groups were defined and investigated in
\cite{luannmbp1,luannmbp2,luddiarijmms03}. But for non-archimedean
fiber bundles they were not yet studied.
\par In this paper wrap groups are defined and studied for fiber
bundles. These fiber bundles are considered over infinite fields
with non-trivial ultra-norms. Then Cayley-Dickson algebras over such
fields are also considered. Fiber bundles over them are introduced.
For such fiber bundles geometric wrap groups are constructed as
well. Their structure is investigated. Such investigation is
motivated by the following reasons. Geometric loops are used in
quantum field theory and they were introduced by Wilson and his
followers in physics, for example, to describe confinement of quarks
\cite{wilson74,fredemar,harvey}.
\par Non-archimedean functional analysis and quantum mechanics develop
intensively in recent times \cite{roo,vla3}. This is stimulated by
several problems. One of them consists in the divergence of some
important integrals and series in the real or complex cases and
their convergence in the non-archimedean case. Therefore, it is
important to consider non-Archimedean wrap semigroups and groups,
that are new objects. There are many principal differences between
classical functional analysis (over the fields $\bf R$ or $\bf C$)
and the non-archimedean one \cite{roo,sch1,wei}. \par The notions of
wrap groups and semigroups in the non-archimedean case are used here
in analogy with the case of manifolds over the real field $\bf R$,
but their meaning is quite different, because non-archimedean
manifolds $M$ modeled on ultra-normed spaces are totally
disconnected with the small inductive dimension $ind(M)=0$ (see \S
6.2 and Chapter 7 in \cite{eng}) and real manifolds are locally
connected with $ind(M)\ge 1$. In the real case loop and wrap groups
$G$ are locally connected for $dim_{\bf R}M<dim_{\bf R}N$, but in
the non-archimedean case they are zero-dimensional with $ind(G)=0$,
where $1\le dim_{\bf R}N\le \infty $ is the dimension of the tangent
Banach space $T_xN$ over $\bf R$ for $x\in N$.
\par In this article wrap groups and semigroups are
considered.  The wrap semigroups of manifolds are quotients of
families of mappings $f$ from one non-archimedean manifold $M$ into
another $N$ with $\lim_{x\to s_0}\bar \Phi ^vf(x)=0$ for $0\le v\le
t$ by the corresponding equivalence relations, where $s_0$ and
$y_0=0$ are marked points in $\bar M$ and $N$ respectively, $M=\bar
M\setminus \{ s_0\}$, ${\bar \Phi }^vf$ are continuous extensions of
the partial difference quotients $\Phi ^vf$. Besides locally compact
manifolds also non-locally compact manifolds $M$ and $N$ are
considered. More generally differentiable spaces modeled on locally
convex non-archimedean spaces are also considered. Then
differentiable fiber bundles with a parallel transport structure are
introduced and for them wrap groups are defined and investigated.
Particularly groups of continuous wraps are also considered. We
consider fibers over different infinite ultra-normed locally compact
and non locally compact fields of zero and positive characteristics.
\par In this article over non-archimedean fields apart from works
of others authors over the fields of $\bf R$ or $\bf C$ we do not
impose an additional condition on an operator of a parallel
transport structure related with tangent vectors. The latter was
used that to bind the parallel transport with the covariant
differentiation on a classical manifold. It was useful for geometric
interpretations and for physical applications. In the present work
we elaborate more general construction without such restriction
using an abstract parallel transport structure taking values from a
structure group of a fiber bundle. This produces wider family of
considered objects of wrap groups. In the future work we plan to
bind this with geometry and physical applications over infinite
fields with non-archimedean multiplicative norms imposing additional
conditions on the parallel transport structure.
\par It is demonstrated that for differentiable spaces wrap groups are commutative.
For fiber bundles with non commutative structure groups wrap groups
are generally non commutative.

\par Semigroups and groups of wraps are investigated in \S 3.
\par The wrap groups are generally non locally compact and have a structure
of differentiable groups modeled on differentiable spaces, which may
be in particular continuous spaces and groups. The main results of
this paper are obtained for the first time and are given in Theorems
3.3, 3.6, 3.11, Proposition 3.13 and Corollaries 3.5, 3.8 and 3.9.

\section{Non-archimedean fiber bundles.}
\par To avoid misunderstandings we first present our definitions
and notations.
\par {\bf 1. Remark.} Let $A$ be an algebra with unit $1$ over a field
$\bf K$ and $\delta $ be some element of $\bf K$ different from
zero. Suppose that $A$ is supplied with a $\bf K$ linear mapping
$x\mapsto x^*$ being an involution such that \par $(1)$ $(x^*)^* =
x$, $x+x^* = tr (x) \in \bf K$, $xx^* = n(x)\in \bf K$. \\ We take
the direct sum of $\bf K$ linear spaces $A_1 := A\oplus A$ and
define the multiplication: \par $(2)$ $(a_1,b_1) (b_1,b_2) = (a_1b_1
- \delta b_2a_2^*, a_1^* b_2 + b_1 a_2)$\\
for each $a_1, a_2, b_1, b_2 \in A$. This multiplication supplies
$A_1$ with the algebraic structure. Certainly, the initial algebra
$A$ is embedded as the subalgebra into $A_1$,
\par $(3)$ $A\ni x\mapsto (x,0)\in A_1$, $(a_1,a_2)^* = (a_1^*, -a_2)$;
\par $(4)$ $tr (a_1, a_2) = tr (a_1)$, $n(a_1a_2) = n(a_1) + \delta n(a_2)$.
\par This doubling procedure applied by induction gives a sequence of embedded
subalgebras ${\bf K} =: {\cal A}_0 \hookrightarrow {\cal A}_1
\hookrightarrow {\cal A}_2\hookrightarrow {\cal A}_3 \hookrightarrow
... $. Henceforward, we consider commutative fields $\bf K$.
Cayley-Dickson algebras ${\cal A}_r$ with $r\le 3$ are alternative:
\par $(5)$ $a(bb) = (ab)b$ and $b(ba) = (bb)a$ for all $a, b\in
{\cal A}_3$, moreover, the quadratic form $n(x)$ is multiplicative
\par $(6)$ $n(ab) = n(a) n(b)$ for all $a, b\in {\cal A}_3$.
\par  In this case the algebra ${\cal A}_1(\alpha )$ is commutative, the algebra
${\cal A}_2(\alpha ,\beta )$ is associative and is called the
algebra of (generalized) quaternions. The Cayley-Dickson algebra
${\cal A}_3(\alpha ,\beta ,\gamma )$ is alternative and simple with
the center $Z({\cal A}_3) = {\bf K}$. But generally the
Cayley-Dickson algebra ${\cal A}_3$ is non-associative.
\par It is possible to take as ${\cal A}_1$ an algebra with a basis $
\{ 1, u \} $ and with the multiplication $u^2=u+\alpha $, where
$\alpha \in \bf K$, $4\alpha +1\ne 0$, and with the involution
$1^*=1$, $u^* = 1-u$.
\par We consider fields $\bf K$ of characteristic $char ({\bf K}) \ne
2$. In this case the Cayley-Dickson algebra ${\cal A}_3$ has a basis
of generators $ \{ 1, u_1,..., u_7 \} $ so that
\par $(7)$ $1^* =1$, $u_j^* = - u_j$ for every $j=1,...,7$, $u_1^2 = - q_1 $,
$u_2^2 = - q_2 $, $u_4^2 = - q_3 $, $u_3^2= - q_1q_2$, $u_5^2 = -
q_1q_3$, $u_6^2 = - q_2q_3$, $u_7^2 = - q_1q_2q_3$, $u_1u_2=u_3$,
$u_1u_4 = - u_5$, $u_2u_4= - u_6$, $u_1u_6 = u_3u_4 = - u_2u_5= -
u_7$, $1u_j= u_j1 = u_j$ and $u_ju_k= - u_ku_j$ for each $1\le j\ne
k$, where $q_1, q_2, q_3 \in \bf K$ (see \cite{dickson,shaefer}). In
more details this algebra is denoted by ${\cal A}_3(q_1, q_2, q_3)$
and its subalgebras are denoted by ${\cal A}_2(q_1, q_2)$, ${\cal
A}_1(q_1)$. The Cayley-Dickson algebra ${\cal A}_3$ has the division
property $ab\ne 0$ for each non zero elements $a, b\in {\cal
A}_3\setminus \{ 0 \} $ if and only if the quadratic form $n(b)$ is
non zero on ${\cal A}_3\setminus \{ 0 \} $.
\par Henceforth, we consider the Cayley-Dickson algebras ${\cal A}_3(q_1, q_2, q_3)$
so that they are alternative with the division property if something
other is not specified. This implies that $G := {\cal A}_3(q_1, q_2,
q_3)\setminus \{ 0 \} $ has the properties $(G1-G4)$:
\par $(G1)$ there is a binary operation $ab \in G$ for all $a, b \in G$;
\par $(G2)$ $1b=b1=b$ for all $b\in G$, that is $1=e$ is the unit
element;
\par $(G3)$ each $a\in G$ has the inverse element $a^{-1}$ so that
$a^{-1}a=aa^{-1} = e$;
\par $(G4)$ $(ab)b = a(bb)$ and $b(ba) = (bb)a$, $(ab)b^{-1} = a$
and $b^{-1}(ba) = a$ for all $a, b \in G$.
\par Property $(G4)$ is called an alternativity, or a weak
associativity. For usual groups the axiom $(G4)$ is replaced on the
associativity: \par $(G5)$ $(ab)c= a(bc)$ for all $a, b , c\in G$.
\\  We shall call $G$ an alternative (or weak associative) group, when
it satisfies Conditions $(G1-G4)$. A usual group satisfies
$(G1-G3,G5)$, so that Condition $(G4)$ follows from $(G5)$. For
short we call a group $G$ in both cases, when a situation is clear.
\par It is necessary to note that an object is known which is called
a path group (a group of paths) at first introduced in physical
literature and then in mathematical. This term does not mean an
algebraic group or a topological group, because compositions are
defined not for all elements, that is only properties $(G2,G3)$ are
fulfilled for the path group while $(G1,G5)$ may be satisfied for
definite combinations of elements only. Wrap groups compose the main
subject of this paper such that they will satisfy $(G1-G5)$ or
$(G1-G4)$.

\par {\bf 2. Definitions.} We consider an infinite field $\bf K$ with a
non trivial non archimedean normalization. We suppose also that $X$
and $Y$ are topological vector spaces over $\bf K$ and $U$ is an
open subset in $X$. For a function $f: U\to Y$ we consider the
associated function
\par $f^{[1]}(x,v,t) := [f(x+tv) - f(x)]/t$ \\
on a set $U^{[1]}$ at first for $t\ne 0$ such that $U^{[1]} := \{
(x,v,t): ~  (x,v,t)\in X^2\times {\bf K}; x\in U, x+tv\in U \} $. If
the function $f$ is continuous on $U$ and $f^{[1]}$ has a continuous
extension on $U^{[1]}$, then we say, that $f$ is continuously
differentiable or belongs to the class $C^1$. The $\bf K$-linear
space of all such continuously differentiable functions $f$ on $U$
is denoted $C^{[1]}(U,Y)$. Then we define by induction functions
$f^{[n+1]}:= (f^{[n]})^{[1]}$ and their spaces $C^{[n+1]}(U,Y)$ for
$n=1,2,3,...$, where $f^{[0]}:=f$, $f^{[n+1]}\in C^{[n+1]}(U,Y)$ has
as the domain $U^{[n+1]} := (U^{[n]})^{[1]}$, $f^{[n]} =: \Upsilon
^n f$.
\par At the same time a differential $df(x): X\to Y$ is defined as
$df(x)v := f^{[1]}(x,v,0)$.
\par Define also partial difference quotient operators $\Phi ^n$
by variables corresponding to $x$ only such that
\par $\Phi ^1f(x;v;t) = f^{[1]}(x,v,t)$ \\
at first for $t\ne 0$ and if $\Phi ^1f$ is continuous for $t\ne 0$
and has a continuous extension on $U^{[1]}=:U^{(1)}$, then we denote
it by ${\bar {\Phi }}^1f(x;v;t)$. We define by induction \par $\Phi
^{n+1} f(x;v_1,...,v_{n+1};t_1,...,t_{n+1}):= \Phi ^1(\Phi
^nf(x;v_1,...,v_n;t_1,...,t_n))(x;v_{n+1};t_{n+1})$ \\
at first for $t_1\ne 0,...,t_{n+1}\ne 0$ on $U^{(n+1)}:= \{
(x;v_1,...,v_{n+1};t_1,...,t_{n+1}): x\in U; v_1,...,v_{n+1}\in X;
t_1,...,t_{n+1}\in {\bf K}; x+v_1t_1\in
U,...,x+v_1t_1+...+v_{n+1}t_{n+1}\in U \} $. If $f$ is continuous on
$U$ and partial difference quotients $\Phi ^1f$,...,$\Phi ^{n+1}f$
have continuous extensions denoted by ${\bar {\Phi }}^1f$,...,
${\bar {\Phi }}^{n+1}f$ on $U^{(1)}$,...,$U^{(n+1)}$ respectively,
then we say that $f$ is of class of smoothness $C^{n+1}$. The $\bf
K$ linear space of all $C^{n+1}$ functions on $U$ is denoted by
$C^{n+1}(U,Y)$, where $\Phi ^0f := f$, $C^0(U,Y)$ is the space of
all continuous functions $f: U\to Y$. \par Then the $n$-th
differential is given by the equation $d^nf(x).(v_1,...,v_n) := n!
{\bar {\Phi }}^nf(x;v_1,...,v_n;0,...,0)$, where $n\ge 1$, also
denote $D^nf=d^nf$. Shortly we shall write the argument of $f^{[n]}$
as $x^{[n]}\in U^{[n]}$ and of ${\bar {\Phi }}^nf$ as $x^{(n)}\in
U^{(n)}$, where $x^{[0]}=x^{(0)}=x$, $x^{[1]}=x^{(1)}=(x,v,t)$,
$v^{[0]}=v^{(0)}=v$, $t_1=t$, $x^{[k]}=(x^{[k-1]},v^{[k-1]},t_k)$
for each $k\ge 1$; $x^{(k)} := (x;v_1,...,v_k;t_1,...,t_k)$.
\par We denote by $C^n_b(U,Y)$ or $C^{[n]}_b(U,Y)$
respectively subspaces of uniformly $C^n$ or $C^{[n]}$ bounded
continuous functions together with ${\bar {\Phi }}^kf$ or $\Upsilon
^kf$ on bounded open subsets of $U$ and $U^{(k)}$ or $U^{[k]}$ for
$k=1,...,n$. These spaces of differentiable functions were
investigated in details in \cite{ludfnafjms09,luseamb}.

\par  We denote by $L(X,Y)$ the space of all
continuous $\bf K$-linear mappings $A: X\to Y$. By $L_n(X^{\otimes
n},Y)$ we denote the space of all continuous $\bf K$ $n$-linear
mappings $A: X^{\otimes n}\to Y$, particularly,
$L(X,Y)=L_1(X^{\otimes 1},Y)$. If $X$ and $Y$ are normed spaces,
then $L_n(X^{\otimes n},Y)$ is supplied with the operator norm: $ \|
A \| := \sup_{h_1\ne 0,...,h_n\ne 0; h_1,...,h_n\in X} \|
A.(h_1,...,h_n) \|_Y/ (\| h_1 \| _X... \| h_n \| _X)$.

\par {\bf 3. Definitions.} Suppose that $M$ is a manifold modeled on a
topological vector space $X$ over $\bf K$ such that its atlas $At
(M) := \{ (U_j,\mbox{ }_M\phi _j): j\in \Lambda _M \} $ is of class
$C^{\alpha ' }_{\beta }$, that is the following four conditions are
satisfied:
\par $(M1)$ $\{ U_j: j\in \Lambda _M \} $ is an open covering of $M$,
$U_j=\mbox{ }_MU_j$; \par $(M2)$ $\bigcup_{j\in \Lambda _M} U_j=M$;
\par $(M3)$ $\mbox{ }_M\phi _j :=\phi _j:
U_j\to \phi _j(U_j)$ is a homeomorphism for each $j\in \Lambda _M$,
$\phi _j(U_j)\subset X$, every $\phi _j(U_j)$ is open in $X$;
\par $(M4)$ $\phi _j\circ \phi _i^{-1}\in C^{\alpha '}_{\beta }$ on
its domain for each $U_i\cap
U_j\ne \emptyset $, \\
where $\Lambda _M$ is a set, $C^{\infty }:= \bigcap_{l=1}^{\infty
}C^l_{\beta }$, $C^{[\infty ]}_{\beta }:=\bigcap_{l=1}^{\infty
}C^{[l]}$, $\alpha '\in \{ n, [n]: 1\le n \le \infty \} $, $\beta
\in \{ \emptyset , b \} $, $C^{\alpha '}_{ \emptyset } := C^{\alpha
'}$.
\par Supply $C^{\alpha }_{\beta }(U,Y)$ with the bounded-open
$C^{\alpha }_{\beta }$ topology denoted by $\tau _{\alpha, \beta }$
generally or $\tau _{\alpha }$ for $\beta =\emptyset $ (or for
compact $U$) with the base \par $(B1)$ $W(P,V)= \{ f\in C^{\alpha
}_{\beta }(X,Y): S^kf|_P\in V, k=0,...,n \} $ \\ of neighborhoods of
zero, where $P$ is bounded and open in $U\subset X$, $P\subset U$,
$V$ is open in $Y$, $0\in V$, $S^k={\bar {\Phi }}^k$ or
$S^k=\Upsilon ^k$ for $\alpha =n$ or $\alpha =[n]$ respectively,
$v_1,...,v_n\in (P-y_0)$, $v^{[k]}_l\in (P-y_0)$ for each $k, l$ for
some marked $y_0\in P$ and $|t_j|\le 1$ for every $j$.
\par If $M$ and $N$ are $C^{\alpha '}_{\beta }$ manifolds
on topological vector spaces $X$ and $Y$ over $\bf K$ respectively,
then let us consider the uniform space $C^{\alpha }_{\beta }(M,N)$
of all mappings $f: M\to N$ such that $f_{j,i}\in C^{\alpha }_{\beta
}$ on its domain for each $j\in \Lambda _N$, $i\in \Lambda _M$,
where $f_{j,i}:= \mbox{ }_N\phi _j\circ f\circ \mbox{ }_M\phi
_i^{-1}$ is with values in $Y$, $\alpha \le \alpha '$. The
uniformity in $C^{\alpha }_{\beta }(M,N)$ is inherited from the
uniformity in $C^{\alpha }_{\beta }(X,Y)$ with the help of charts of
atlases of $M$ and $N$. If $M$ is compact, then $C^{\alpha }_b(M,N)$
and $C^{\alpha }(M,N)$ coincide.
\par The family of all homeomorphisms $f: M\to M$ of class
$C^{\alpha }_{\beta }$ is denoted by $Diff^{\alpha }_{\beta }(M)$.
\par Let $\gamma $ be a set, then we denote by $c_0(\gamma ,{\bf K})$
the normed space consisting of all vectors $ x= \{ x_j\in {\bf K}:
j\in \gamma , \mbox{ for each } \epsilon >0 \mbox{ the set } \{ j:
|x_j|>\epsilon \} \mbox{ is finite } \} $, where \par $(N1)$ $\| x
\| := \sup_{j\in \gamma } |x_j|$. \\ In view of the Kuratowski-Zorn
lemma it is convenient to consider $\gamma $ as an ordinal.
Henceforth, suppose that $X=c_0(\gamma _X,{\bf K})$ and
$Y=c_0(\gamma _Y,{\bf K})$.
\par As a generalization let $X=c_0(\gamma ,{\cal A}_r)$ be a normed
space over ${\cal A}_r = {\cal A}_r(q_1,...,q_r)$ consisting of all
vectors of the form $ x= \{ x_j\in {\cal A}_r: j\in \gamma , \mbox{
for each } \epsilon >0 \mbox{ the set } \{ j: |x_j|>\epsilon \}
\mbox{ is finite } \} $, where \par $(N2)$ $\| x \| := \sup_{j\in
\gamma } |x_j|_r$, where $1\le r \le 3$ (see \S 1),
\par $(N3)$ $|x_j| := |x_j|_r := \max_{k=0,...,2^r-1} |\mbox{}_kx_j|$, $x_j =
\mbox{}_0x_ju_0+...+\mbox{}_{2^r-1}x_j u_{2^r-1}$, $\mbox{}_kx_j\in
\bf K$ for all $k, j$, $u_0=1$, $ \{ u_0, ...,u_{2^r-1} \} $ is a
basis of generators of ${\cal A}_r$.
\par By a vector space $X$ over ${\cal A}_r$ we shall undermine the direct sum $X=
\mbox{}_0X_0u_0\oplus ... \oplus \mbox{}_{2^r-1}X u_{2^r-1}$, where
$\mbox{}_0X,...,\mbox{}_{2^r-1}X$ are pairwise isomorphic linear
spaces over a field $\bf K$, $1\le r\le 3$. That is the following
conditions $(L1-L4)$ are satisfied:
\par $(L1)$ $X$ is the additive commutative group,
\par $(L2)$ $a(x+y)=ax+ay$ and $(x+y)a=xa+ya$,
\par $(L3)$ $(a+b)x = ax+bx$ and $x(a+b)=xa+xb$ for all $a, b\in {\cal
A}_r$ and $x, y \in X$,
\par $(L4)$ $(ab)\mbox{}_0x = a(b\mbox{}_0x)$ and $\mbox{}_0x (ab) =
(\mbox{}_0xa)b$ for all $a, b\in {\cal A}_r$ and $\mbox{}_0x\in
\mbox{}_0X$. \par We consider a vector space $X$ over ${\cal A}_r$
supplied with a topology $\tau $ with jointly continuous operations
of addition of vectors $X^2\ni (x,y)\mapsto x+y\in X$ and their
multiplication ${\cal A}_r\times X\ni (a,x)\mapsto ax\in X$ and
${\cal A}_r\times X\ni (a,x)\mapsto xa\in X$ on scalars from ${\cal
A}_r$ relative to $\tau $ and the norm topology $|a|=|a|_r$ in
${\cal A}_r$. Such $X$ is called a topological vector space.
\par The algebras ${\cal A}_1$ and ${\cal A}_2$ are associative
and for them Conditions $(L1,L3,L4)$ imply that $(ab)x=a(bx)$ and
$x(ab)=(xa)b$ for all $x\in X$ and $a, b$ in ${\cal A}_1$ or ${\cal
A}_2$ respectively, when $1\le r \le 2$. Moreover, over ${\cal A}_1$
the multiplication on scalars is commutative so that we can consider
the equality $ax=xa$ for all $a\in {\cal A}_1$ and $x\in X$, when
$r=1$. Thus each vector space over ${\cal A}_r$ with $1\le r \le 3$
is also a vector space over the initial field ${\bf K} = {\cal
A}_0$. \par We mention also that algebras ${\cal A}_r$ are finite
dimensional over $\bf K$ and additions and multiplications and
taking of the conjugate $a\mapsto a^*$ in them are continuous
relative to the norm $|a|=|a|_r$. If $r\le 3$, then by our
convention of \S 1 $n(a)=aa^* = a^*a \in \bf K$ and $a^{-1} =
a^*/n(a)$ for $a\ne 0$, since $\alpha ^* = \alpha $ for each $\alpha
\in \bf K$. This implies that for the considered Cayley-Dickson
algebras the inversion ${\cal A}_3\setminus \{ 0 \} \ni a\mapsto
a^{-1}\in {\cal A}_3\setminus \{ 0 \}$ is also continuous.
\par The family of all continuous ${\cal A}_r-$additive $\bf K$-linear
mappings $A: X\to Y$ for topological vector spaces $X$ and $Y$ over
${\cal A}_r$ we denote by $K_q(X,Y)$, $1\le r \le 3$. Their
subfamily of right ${\cal A}_r-$linear mappings $A(\mbox{}_0xb) =
(A\mbox{}_0x)b$ for all $\mbox{}_0x\in \mbox{}_0X$ and $b\in {\cal
A}_r$ we denote by $K_r(X,Y)$ or $L(X,Y)$. The subfamily of all left
${\cal A}_r-$linear mappings $A(b\mbox{}_0x) = b(A\mbox{}_0x)$ for
all $\mbox{}_0x \in \mbox{}_0X$ and $b\in {\cal A}_r$ will be
denoted by $K_l(X,Y)$. Certainly that over ${\cal A}_1$ these spaces
coincide $K_r(X,Y)=K_l(X,Y) =K_q(X,Y)$, since the algebra ${\cal
A}_1$ is commutative and associative.
\par Each topological vector space $X$ over ${\cal A}_r$ with $1\le r \le 3$ is a
topological vector space $X_{\bf K}$ over a field $\bf K$ as well.
So by a manifold $M$ modeled on $X$ we shall undermine the
corresponding manifold $M_{\bf K}$ modeled on $X_{\bf K}$ such that
a class of smoothness $C^{\alpha }_{\beta }$ of $M$ is that of
$M_{\bf K}$, but with one additional condition:
\par $(M5)$ each differential $d(\phi _j\circ \phi _i^{-1}(x))$ is a right
${\cal A}_r-$linear operator, that is, it belongs to $L(X,X)$ for
each $x$ in its domain for each $U_i\cap U_j\ne \emptyset $.
\par Therefore we obtain also as in \S 3 above uniform spaces $C^{\alpha }_{\beta }(M,N)$
for manifolds $M$ and $N$ modeled on topological vector spaces $X$
and $Y$ over ${\cal A}_r$ with $1\le r \le 3$. Certainly the
uniformity in $C^{\alpha }_{\beta }(M,N)$ is inherited from the
uniformity in $C^{\alpha }_{\beta }(X,Y)$ with the help of charts of
atlases of $M$ and $N$.

\par {\bf 4. Plots.} Since the Cayley-Dickson algebra ${\cal A}_3$
is non-associative, we consider a non-associative subgroup $G$ of
the family $Mat_q({\bf O})$ of all square $q\times q$ matrices with
entries in ${\cal A}_3$. More generally $G$ is a group which has a
$C^{\alpha }_{\beta }$ manifold  structure over ${\cal A}_r$ and
group's operations are $C^{\alpha }_{\beta }$ mappings (see also \S
1), where $0\le r \le 3$. Such groups $G$ we call also a $C^{\alpha
}_{\beta }$ Lie group over ${\cal A}_r$. In particular, it may be a
multiplicative or and additive subgroup of $L(X,X)$ (see \S 3).

\par As a generalization of manifolds we use the
following (over $\bf R$ and $\bf C$ see \cite{gajer,souriau}). We
adopt that a subset $C$ of a vector space $X$ over ${\cal A}_r$ is
called ${\cal A}_r$ absolutely convex if $ax+by$ and $xa+yb\in C$
for all $x, y \in C$ and $a, b \in {\cal A}_r$ with $|a|\le 1$ and
$|b|\le 1$. Translates $(z+C)$ of absolutely ${\cal A}_r$ convex
sets $C$ are called ${\cal A}_r$ convex, where $z\in X$. A
topological vector space $X$ over ${\cal A}_r$ is called locally
${\cal A}_r$ convex if it has a base of ${\cal A}_r$ convex
neighborhoods of zero.
\par Let $X$, $Y$ be topological vector spaces over ${\cal A}_r$, $0\le
r \le 3$, particularly ${\cal A}_0= {\bf K} $. Suppose that $M$ is a
Hausdorff topological space supplied with a family $ {\cal P}_M :=
\{ h : U\to M \} $ of the so called plots $h$ which are continuous
maps satisfying conditions $(D1-D5)$:
\par $(D1)$ each plot has as a domain an ${\cal A}_r$ convex subset $U= \mbox{}_hU$ in $X$;
\par $(D2)$ if $h: U\to M$ is a plot, $V$ is an ${\cal A}_r$ convex subset in
$Y$ and $g: V\to U$ is an $C^{\alpha }_{\beta }$ mapping, then
$h\circ g$ is also a plot;
\par $(D3)$ every constant map from an ${\cal A}_r$ convex set $U$ in $X$ into $M$ is a plot;
\par $(D4)$ if $U$ is an ${\cal A}_r$ convex set in $X$ and
$ \{ U_j: j\in J \} $ is a covering of $U$ by ${\cal A}_r$ convex
sets in $X$, each $U_j$ is open in $U$, $h: U\to M$ is such that
each its restriction $h|_{U_j}$ is a plot, then $h$ is a plot. We
suppose of course that \par $(D5)$ the family of subsets $h_k(U_k)$
which are ranges of plots, $k \in \lambda _M$ for all plots $h_k$ of
$M$ forms a covering of $M$, but not necessarily open, where
$\lambda _M$ is a set. \par Then $M$ is called a $C^{\alpha }_{\beta
}$-differentiable space.
\par A mapping $f: M\to N$ between two $C^{\alpha }_{\beta }$-differentiable spaces
$M$ and $N$ is called $C^{\alpha }_{\beta }$ differentiable if it
continuous and for each plot $h: U\to M$ the composition $f\circ h:
U\to N$ is a $C^{\alpha }_{\beta }$ plot of $N$.
\par To supply a family of $C^{\alpha }_{\beta }(M,N)$ mappings between
$C^{\alpha }_{\beta }$ differentiable spaces with a uniformity we
use the following particular case. We suppose that $M$ and $N$ have
families of plots ${\cal T}_M = \{ h_j^M\in {\cal P}_M: j\in \Lambda
_M \} $ and ${\cal T}_N = \{ h_j^N\in {\cal P}_N: j\in \Lambda _N \}
$ correspondingly, so that domains of plots $h_j^M(U^M_j)$ and
$h_k^N(U^N_k)$ for $M$ with $j\in \Lambda _M$ and $N$ with $k\in
\Lambda _N$ respectively form coverings of $M$ and $N$ satisfying
Conditions $(D6-D8)$:
\par $(D6)$ $h_j^M: U^M_j\to h_j^M(U^M_j)$ and $h_k^N: U^N_k\to h^N_k(U^N_k)$ are
bijective for all $j\in \Lambda _M$ and $k\in \Lambda _N$, where
$\Lambda _M\subset \lambda _M$ and $\Lambda _N\subset \lambda _N$
are subsets;
\par $(D7)$ transition mappings $(h_l^M)^{-1}\circ h_j^M$ and
$(h_s^N)^{-1}\circ h_k^N$ are $C^{\alpha }_{\beta }$ differentiable
mappings on their domains for each $h_l^M(U^M_l)\cap h_j^M(U^M_j)\ne
\emptyset $ and $h_s^N(U^N_s)\cap h_k^N(U^N_k)\ne \emptyset $ with
$l, j \in \Lambda _M$, $s, k \in \Lambda _N$;
\par $(D8)$ supply $M$ and $N$ with topologies $\tau _{p,M}$ and
$\tau _{p,N}$ having bases $h_j^M(V)$ for each $V$ open in $U^M_j$
and every $j\in \Lambda _M$ for $M$, while $h_k^N(V)$ for each $V$
open in $U^N_k$ and every $k\in \Lambda _N$ for $N$ correspondingly.
\par Then we form atlases for $(M, \tau _{p,M})$ and $(N, \tau
_{p,N})$ with (generalized) charts $(h_j^M(U_j^M), (h_j^M)^{-1})$
and $(h_k^N(U^N_k), (h_k^N)^{-1})$ with $j\in \Lambda _M$ and $k\in
\Lambda _N$ respectively. \par Thus we get the uniform space
$C^{\alpha }_{\beta }(M,N)$ using transition mappings between
(generalized) charts as in \S 3. That is $C^{\alpha }_{\beta }(M,N)$
consists of all mappings $f: M\to N$ such that $f_{k,j}\in C^{\alpha
}_{\beta }$ on its domain for each $k\in \Lambda _N$, $j\in \Lambda
_M$, where $f_{k,j} := (h^N_k)^{-1}\circ f\circ h^M_j$ is with
values in $Y$, $\alpha \le \alpha '$. Here $U^M_j$ is considered as
open in the vector space $X_j := span_{{\cal A}_r} U^M_j$. For this
we supply $X_j$ with a base of topology generated by neighborhoods
of zero $\lambda (W-x)$ for each $x\in U^M_j$ so that $U^M_j$ is
absolutely convex in $X$, $W$ is open and bounded in $U^M_j$
relative to the topology in $U^M_j$ inherited from $X$, $\lambda \in
{\cal A}_r$ with $0< |\lambda |<1$. The uniformity in $C^{\alpha
}_{\beta }(M,N)$ is inherited from the uniformity in $C^{\alpha
}_{\beta }(X,Y)$ with the help of (generalized) charts of atlases of
$M$ and $N$. This means that a base of entourages of the diagonal
(see also Chapter 8 in \cite{eng}) in $C^{\alpha }_{\beta }(M,N)$ is
formed by sets $ \{ f, g \in C^{\alpha }_{\beta }(M,N): (f_{l,j} -
g_{l,j}) \in W(P,V) ~ \forall l \} $, where $P$ is a bounded open
subset in $U^M_j$ and $V$ is an open neighborhood of zero in $Y$,
$j\in \Lambda _M$ and $l\in \Lambda _N$. In $C^{\alpha }_{\beta
}(U,Y)$ the base $W(P,V)$ was defined in $(B1)$ \S 3. \par  A
topological group $G$ is called an $C^{\alpha }_{\beta
}$-differentiable group if its group operations are $C^{\alpha
}_{\beta }$-differentiable mappings.

\par {\bf 5. Transformation groups.} A pair $(G,F)$  is called a $C^{\alpha }_{\beta
}$ transformation group acting from the left if the following four
conditions are satisfied:
\par $(T1)$ $G$ is a $C^{\alpha }_{\beta }$ differentiable group over ${\cal A}_r$,
where a number $r$ is supplied, $0\le r \le 3$;
\par $(T2)$ $F$ is a  $C^{\alpha }_{\beta }$ differentiable space over ${\cal
A}_r$;
\par $(T3)$ a $C^{\alpha }_{\beta }$ mapping $(g,x)\in G\times F\to gx \in F$ is
given;
\par $(T4)$ for each $g\in G$ a $C^{\alpha }_{\beta }$ mapping $l_g: F\to F$
is defined so that $l_g(x)=gx$ and $l_g\in Diff^{\alpha }_{\beta
}(F)$, $l_{gh} = l_g\circ l_h$.
\par If $l_{gh} = l_h\circ l_g$, then $G$ acts on $F$ from the
right. In particular, if $G$ is a $C^{\alpha }_{\beta }$ Lie group
and $N$ is a $C^{\alpha }_{\beta }$ manifold, then $(G,F)$ is called
a $C^{\alpha }_{\beta }$ transformation Lie group.
\par A transformation group $(G,F)$ with $l_g = id_F$ for each $g\in
G$ is called trivial, where $id_F(y) =y$ for each $y\in F$.
Contrary, when $l_g=id_F$ if and only if $g$ is a unit element of
$G$, $g=e$, the transformation group $G$ is called effective.
\par {\bf 6. Fiber bundles.}
\par  Let $E$, $N$, $F$ be all either $C^{\alpha '}_{\beta }$-manifolds or
$C^{\alpha '}_{\beta }$-differentiable spaces over ${\cal A}_r$ with
$0\le r \le 3$. Certainly a manifold is a particular case of a
differentiable space. Let also $G$ be a $C^{\alpha '}_{\beta }$
group over ${\cal A}_r$, $\alpha \le \alpha '\le \infty $. We
suppose that a projection $\pi : E\to N$ is given together with an
atlas $\Psi = \{ \psi _j \} $ of $E$ so that
\par $(F1)$ to each chart $\psi _j \in \Psi $ an open subset $V_j$
in $N$ is counterposed and \par $(F2)$ the mapping $\psi _j: \pi
^{-1}(V_j)\to V_j\times F$ is the $C^{\alpha '}_{\beta }$
diffeomorphism so that $\psi _j (\pi ^{-1}(x)) = \{ x \} \times F$:
\par $pr_{V_j}\circ \psi _j = \pi |_{\pi ^{-1}(V_j)}$, \\ where
$pr_{V_j}(x\times y) =x$ for each $x\in V_j$ and $y\in F$;
\par $(F3)$ a system of open subsets $ \{ V_j: j\in J \} $ forms a
covering of $N$. We get from $(F1-F3)$ that $\pi : E\to N$ is open
and surjective. Moreover,
\par $(F4)$ the mapping $\psi _{j,x} = pr_F\circ \psi _j |_{\pi
^{-1}(x)} : \pi ^{-1}(x)\to F$ defines the $C^{\alpha '}_{\beta }$
diffeomorphism of the fiber $F_x := \pi ^{-1}(x)$ on the typical
fiber $F$, where $pr_F: (x\times y)=y$ for all $x\in V_j$ and $y\in
F$.
\par Using restrictions of mappings $\psi _j$ we can choose $V_j$ as
domains in $N$ either $h_j(U_j)$ of plots from ${\cal T}_N$ (see \S
4) in the case of the differentiable space or of charts in the case
of the manifold. Thus either plots or charts on $N\times F$ are
transferred by $\psi _j^{-1}$ onto plots or charts respectively on
$\pi ^{-1}(V_j)$.
\par Let $\psi _j, \psi _l\in \Psi $ and $V_j\cap V_k\ne \emptyset
$. In view of $(F4)$ \par $(F5)$ for each $x\in V_j\cap V_k$ the
mapping is defined: $g_{j,k}: V_j\cap V_k \ni x \mapsto g_{j,k}(x) =
\psi _{k,x}\circ \psi _{j,x} ^{-1} \in Diff^{\alpha '}_{\beta }(F)$.
That is to each point $x\in V_j\cap V_k$ a diffeomorphism of $F$
corresponds. Moreover, these mappings $g_{j,k}$ satisfy the
following conditions:
\par $(F6)$ $g_{j,k}(x) = (g_{k,j}(x))^{-1}$, $g_{j,j}(x) = id_F$,
where $id_F(y) =y$ for each $y\in F$, $g_{l,j} (x) = g_{l,k}(x)\circ
g_{k,j}(x)$ for each $x\in V_l\cap V_j\cap V_k$;
\par $(F7)$ a $C^{\alpha '}_{\beta }$ differentiable transformation group $(G,F)$ is
given so that $g_{j,k}(x)\in G$ for each $x\in V_j\cap V_k$, when
$V_k\cap V_j\ne \emptyset $, $g_{j,j}(x) = e\in G$, where $e$
denotes the unit element in $G$, also $l_{g_{j,k}(x)} = \psi
_{j,x}\circ \psi _{k,x}^{-1}$ (see \S 5 as well).
\par If Conditions $(F1-F7)$ are satisfied then $E(N,F,G,\pi ,\Psi )$
is called a fiber bundle with a fiber space $E$, a base space $N$, a
typical fiber $F$, projection $\pi $ and a structural group $G$ over
${\cal A}_r$, and an atlas $\Psi $, while the mappings $g_{j,k}$ are
called the transition functions.

\par Local trivializations $\phi _j\circ \pi \circ \Psi _k^{-1}:
V_k(E)\to V_j(N)$ induce the $C^{\alpha '}_{\beta }$-uniformity in
the family $\cal W$ of all principal $C^{\alpha '}_{\beta }$-fiber
bundles $E(N,F,G,\pi ,\Psi )$, where $V_k(E) = \Psi
_k(U_k(E))\subset X(N)\times X(F)$, $V_j(N) = \phi _j(U_j(N))\subset
X(N)$, where $X(F)$, $X(G)$ and $X(N)$ are ${\cal A}_r$-vector
spaces on which $F$, $G$ and $N$ are modeled, $(U_k(E),\Psi _k)$ and
$(U_j(N),\phi _j)$ are either ranges of plots or charts of atlases
of $E$ and $N$, $\Psi _k =\Psi _k^{E}$, $\phi _j = \phi _j^{N}$ (see
also \S 4).

\par If $G=F$ and $G$ acts on itself by left shifts, then a fiber bundle is
called the principal fiber bundle and is denoted by $E(N,G,\pi ,\Psi
)$. As an example a multiplicative group $G={\cal A}_r^*$ may be,
where ${\cal A}_r^*$ denotes the multiplicative group ${\cal
A}_r\setminus \{ 0 \} $. If $G=F= \{ e \} $, then $E$ reduces to
$N$.

\section{Wrap groups and semigroups.}
\par {\bf 1. Parallel transport structure.}  Let $cl_X (A)$ denote the closure of a
subset $A$ in $X$. Let $\hat M$ and $\bar M$ be two $C^{\alpha
'}_{\beta }$ differentiable spaces over ${\cal A}_r$ (see Conditions
2.4$(D1-D8)$ above), $0\le r \le 3$, and with marked points $\{
{\hat s}_{0,q}\in {\hat M}_f: q=1,...,2k \} $ and a $C^{\alpha
'}_{\beta }$-mapping $\Xi : {\hat M}\to \bar M$, where ${\hat M}_f$
is a subset defined below, such that the following conditions
$(S1-S5)$ are satisfied:
\par $(S1)$ $\Xi ^{-1}(x)$ consists of one or at most finite number of distinct points for
each $x\in \bar M$, we denote by ${\hat M}_f$ the set of all $y\in
{\hat M}$ so that $y\in \Xi ^{-1}(x)$  for some $x\in {\bar M}$ with
$\Xi ^{-1}(x)$ consisting of a finite number of distinct points more
than one;
\par $(S2)$ $\Xi $ is
surjective and bijective from ${\hat M}\setminus {\hat M}_f$ onto
${\bar M}\setminus {\bar M}_f$ open in $\bar M$, $\Xi ({\hat
s}_{0,q})=\Xi ({\hat s}_{0,k+q}) = s_{0,q}$ for each $q=1,...,k$,
where ${\bar M}_f := \Xi ({\hat M}_f)$;
\par $(S3)$ for each point $x\in \bar M$ or $y\in \hat M$ there
exist ranges of plots $h_{\bar M}(V)$ and $h_{\hat M}(U)$ being open
neighborhoods of $x$ and $y$ respectively in $\bar M$ and $\hat M$
so that $V$ and $U$ are ${\cal A}_r$ convex in a topological vector
space $X$ over ${\cal A}_r$ on which $\bar M$ and $\hat M$ are
modeled so that $h_{\bar M} : V\to h_{\bar M}(V)$ and $h_{\hat M}: U
\to h_{\hat M}(U)$ are bijective;
\par $(S4)$ the closure $cl_{\bar M} \breve{M} = \bar M$ of $\breve{M}$ in $\bar M$
is the entire $\bar M$, where $\breve{M} := \{ x\in {\bar M}: x\in
h_{\bar M}(V)$ $\mbox{for some plot}$ $h\in {\cal T}_{\bar M}$
$\mbox{of}$ $\bar M$ $\mbox{with}$ $V$ $\mbox{open in}$ $X \} $;
\par $(S5)$ $\breve{M} \subset {\bar M}\setminus {\bar M}_f$.
\par By $span_{{\cal A}_r} V$ we denote a vector space consisting of
all finite ${\cal A}_r$ vector combinations of vectors from $V$ and
with multiplication on constants from ${\cal A}_r$, where the
multiplications may be on both sides, when $2\le r \le 3$.
\par Mention that ${\cal A}_r$ with $r\le 3$ are division algebras and
for matrices with entries in ${\cal A}_r$ the Gauss' algorithm is
accomplished, so they have ranks by rows and columns which coincide.
This means that an ${\cal A}_r$ vector independence and a dimension
over ${\cal A}_r$ are well-defined as it was outlined already by
Dickson \cite{dickson}.
\par  Put $M := {\bar M} \setminus \{ s_{0,q}: q=1,...,k \} $.
Here particularly $C^{\alpha }_{\beta }$ manifolds $\hat M$ and
$\bar M$ may be as well.
\par A parallel transport structure on a $C^{\alpha '}_{\beta }$-differentiable
principal $G$-bundle $E(N,G,\pi ,\Psi )$ for $C^{\alpha '}_{\beta }$
differentiable spaces $\bar M$ and $\hat M$ as above over the same
field or an algebra ${\cal A}_r$, $0\le r \le 3$, with $\alpha '\ge
\alpha $ assigns to each $C^{\alpha '}_{\beta }$ mapping $\gamma $
from $\bar M$ into $N$ and points $u_1,...,u_k\in E_{y_0}$, where
$y_0$ is a marked point in $N$, $y_0=\gamma (s_{0,q})$, $q=1,...,k$,
a unique $C^{\alpha }_{\beta }$ mapping ${\bf P}_{{\hat \gamma },u}:
{\hat M}\to E$ satisfying conditions $(P1-P4)$:
\par $(P1)$ take
${\hat \gamma }: {\hat M}\to N$ such that ${\hat \gamma }=\gamma
\circ \Xi $, then ${\bf P}_{{\hat \gamma },u}({\hat s}_{0,q})=u_q$
for each $q=1,...,k$ and $\pi \circ {\bf P}_{{\hat \gamma },u}={\hat
\gamma }$
\par $(P2)$ ${\bf P}_{{\hat \gamma },u}$ is the $C^{\alpha }_{\beta }$-mapping
by $\gamma $ and $u$;
\par $(P3)$ for each $x\in \hat M$ and every $\phi \in Diff^{\alpha }_{\beta }
({\hat M}, \{ {\hat s}_{0,1},...,{\hat s}_{0,2k} \} )$ the equality
${\bf P}_{{\hat \gamma },u}(\phi (x)) = {\bf P}_{{\hat \gamma }\circ
\phi ,u}(x)$ is satisfied, where $Diff^{\alpha }_{\beta }({\hat M},
\{ {\hat s}_{0,1},...,{\hat s}_{0,2k} \} )$ denotes the group of all
$C^{\alpha }_{\beta }$ homeomorphisms of $\hat M$ preserving marked
points $\phi ({\hat s}_{0,q})={\hat s}_{0,q}$ for each $q=1,...,2k$;
\par $(P4)$ ${\bf P}_{{\hat \gamma },u}$ is $G$-equivariant,
which means that ${\bf P}_{{\hat \gamma },uz}(x) = {\bf P}_{{\hat
\gamma },u}(x)z$ for every $x\in {\hat M}$ and each $z\in G$.

\par Two $C^{\alpha '}_{\beta }$-differentiable principal $G$-bundles
$E_1$ and $E_2$ with parallel transport structures $(E_1,{\bf P}_1)$
and $(E_2,{\bf P}_2)$ are called isomorphic, if there exists an
isomorphism $h: E_1\to E_2$ such that ${\bf P}_{2,{\hat \gamma
},u}(x) = h ({\bf P}_{1,{\hat \gamma }, h^{-1}(u)}(x))$ for each
$C^{\alpha }_{\beta }$-mapping $\gamma : {\bar M} \to N$ and $u_q\in
(E_2)_{y_0}$, where $q=1,...,k$,
$h^{-1}(u)=(h^{-1}(u_1),...,h^{-1}(u_k))$.

\par {\bf 2. Subspaces.} For $M = {\bar M}\setminus \{ s_{0,q}: q=1,...,k \} $
either ranges $h_j(W_j)=U_j$ of plots for a $C^{\alpha '}_{\beta }$
differentiable space $\bar M$ or an atlas $At({\bar M})$ of a
$C^{\alpha '}_{\beta }$ manifold $\bar M$ with charts $(U_j,\phi
_j)$, $j\in \Lambda _M$ we put
$$(1)\mbox{ }U_l = {\bar U}_l \setminus \{ s_{0,q}: q=1,...,k \} $$
for each $l=1,...,k$ so that $s_{0,q}\in {\bar U}_q$ for every
$q=1,...,k$ and either
\par $(2)$ $h_q = {\bar h}_q |_{W_q}$ with $W_q = {\bar W}_q\setminus
h_q^{-1}(s_{0,q})$ or $\phi _l = {\bar \phi }_l|_{U_l}$ for all $q,
l =1,...,k$; \par $(3)$ $\mbox{ }U_j={\bar U}_j\mbox{ and } \phi
_j={\bar \phi }_j\mbox{ for each }j>k,$ $\{ s_{0,q}: q=1,...,k \}
\cap {\bar U}_j =\emptyset $ for each $j>k$; $j\in \Lambda _M =
\Lambda _{\bar M}$, where due to the Kuratowski-Zorn theorem
\cite{eng} we can consider, that $\Lambda _M$ is an ordinal.
\par  Let the spaces be the same as in \S 2.3, 4 with the covering
of $M$ defined by Conditions $(1-3)$. Suppose that $M$ is modeled on
$X$ and $N$ on $Y$, where $X$ and $Y$ are ${\cal A}_r$ vector
spaces, $0\le r \le 3$. Then we consider their subspaces of mappings
preserving marked points relative to a given mapping $\theta \in
C^{\alpha }_{\beta }(M,N)$:
$$(4)\mbox{ }C^{\alpha , \theta }_{\beta ,0}((M, \{ s_{0_q}: q=1,..,k \} ); (N,y_0)) :=
\{ f\in C^{\alpha }_{\beta }({\bar M}, N):$$ $$
\lim_{|t_1|+...+|t_m|\to 0}S^m(f_{l,j} - \theta _{l,j})(w^m_q
(t_1,...,t_m)) =0$$ $$ \mbox{ for each } m\in \{ 0,1,...,n \}, ~
\forall j\in \Lambda _M, ~ \forall l \in \Lambda _N, \quad \forall
q=1,...,k \} ,$$ where either $S^m={\bar {\Phi }}^m$ or
$S^m=\Upsilon ^m$ for $\alpha =n$ or $\alpha =[n]$ respectively, an
argument is either $w^m_q(t_1,...,t_m) = x^{(m)}_q\in U^{(m)}_{l,j}$
for $S^m = {\bar {\Phi }}^m$ or $w^m_q (t_1,...,t_m)= x^{[m]}_q\in
U^{[m]}_{l,j}$ for $S^m=\Upsilon ^m$, $f^{[m]} = \Upsilon ^mf$,
where $t_1,...,t_m\in \bf K$, $x^{[0]}_q=x^{(0)}_q=x_q=s_{0,q}$,
$x^{[1]}_q=x^{(1)}_q=(x_q,v,t)$, $v^{[0]}=v^{(0)}=v$, $t_1=t$,
$x^{[m]}_q=(x^{[m-1]}_q,v^{[m-1]},t_m)$ for each $m\ge 1$, $x^{(m)}
:= (x_q;v_1,...,v_m;t_1,...,t_m)$ so that $w^m_q$ is in a domain
either $U^{(m)}_{l,j}$ of ${\bar {\Phi }}^m(f_{l,j} - \theta
_{l,j})$ or $U^{[m]}_{l,j}$ of $\Upsilon ^m(f_{l,j} - \theta
_{l,j})$ correspondingly (see also \S \S 2.2,3). For $\alpha =
\infty $ Condition $(4)$ is imposed for each natural value of $m$.
When points $s_{0,q}$ and $y_0$ are specified we can write shortly
$C^{\alpha , \theta }_{\beta ,0}(M, N)$ instead of $C^{\alpha ,
\theta }_{\beta ,0}((M, \{ s_{0_q}: q=1,..,k \} ); (N,y_0))$.
\par As usually a
diffeomorphism group $Diff^{\alpha }_{\beta }({\bar M})$ of the
differentiable space ${\bar M}$ consists of all surjective bijective
mappings $f$ from $\bar M$ onto $\bar M$ with $f$ and $f^{-1}$
belonging to the differentiability class $C^{\alpha }_{\beta }$. We
consider the following subgroup also:
$$(5)\mbox{ }Diff^{\alpha }_{\beta ,0}(M) := \{ f\in Diff ^{\alpha }_{\beta
}({\bar M}):  f(s_{0,q})=s_{0,q} ~ \forall q=1,...,k \} .$$  ~ ~ ~
We introduce also \par $(6)$ the family $Di^{\alpha }_{\beta ,0}(M)$
of all continuous mappings $f$ from $\bar M$ into $\bar M$ so that
the restriction $f|_{\breve{M}}$ on $\breve{M}$ is bijective and
surjective from $\breve{M}$ onto $\breve{M}$ so that $f$ and
$f^{-1}$ belong to the class $C^{\alpha }_{\beta }$ and
$f(s_{0,q})=s_{0,q}$ for each $q=1,...,k$ (see also \S 1,
$(S4,S5)$). That is if $f\in Di^{\alpha }_{\beta ,0}(M)$, then
$f|_{\breve{M}} \in Diff^{\alpha }_{\beta }(\breve{M})$. We call
such $f$ a (generalized) diffeomorphism of a $C^{\alpha '}_{\beta }$
differentiable space $\bar M$, while $Di^{\alpha }_{\beta ,0}(M)$ we
call a group of (generalized) diffeomorphisms of a differentiable
space $\bar M$ preserving marked points $s_{0,q}$, $q=1,...,k$.
\par It is worth to mention that in the particular case when $\bar
M$ is a $C^{\alpha '}_{\beta }$ manifold with $\alpha \le \alpha '$
the family $Di^{\alpha }_{\beta ,0}(M)$ coincides with $Diff^{\alpha
}_{\beta ,0}(M)$, since for the manifold $\bar M$ all charts $\phi
_j(U_j)$ are open in $X$, so the set $\breve{M}$ is the entire $\bar
M$ (see also \S 2.3). \par  The action of $Di^{\alpha }_{\beta
,0}(M)$ on $M$ induces isomorphism classes of $C^{\alpha }_{\beta }$
principal $G$ fiber bundles with parallel transport structure for
which mappings $\gamma : {\bar M}\to N$ belong to the uniform space
$C^{\alpha , w_0}_{\beta ,0}((M, \{ s_{0_q}: q=1,..,k \} );
(N,y_0))$, where a mapping $\theta = w_0$ is constant: $w_0({\bar
M})= \{ y_0 \} $ (see \S 1 and \S 2.6). We denote by $(S^ME)_{\alpha
, \beta } := (S^{M, \{ s_{0,q}: q=1,...,k \} }E; (N,y_0),G,{\bf
P})_{\alpha , \beta }$ a set of $C^{\alpha }_{\beta }$-closures of
all such isomorphism classes.
\par Recall that a subset $A$ of a topological space $B$ so that $A$ is dense
in itself and closed in $B$ is called a perfect set $A$ \cite{eng}.
\par A topological space $X$ is called a $T_0$ space if for each $x\ne y\in X$
there exists an open subset containing only one of these two points.
A topological space $X$ is called a $T_1$-space if for each pair of
distinct points $x\ne y \in X$ an open subset $U$ of $X$ exists so
that $x\in U$ and $y\notin U$.

\par {\bf 3. Theorems.} {\it {\bf 1.} A uniform space
$(S^ME)_{\alpha ,\beta }$ from \S 2 exists and it has a structure of
a topological $T_1$ alternative monoid with a unit and with a
cancelation property and a multiplication operation of $C^l_{\beta
}$ class with $l=\alpha ' - \alpha $ ($l=\infty $ for $\alpha
'=\infty $). If $M$, $N$ and $G$ are separable, then $(S^ME)_{\alpha
,\beta }$ is separable. If $N$ and $G$ are complete, then
$(S^ME)_{\alpha ,\beta }$ is complete.
\par {\bf 2.} If $G$ is associative, then $(S^ME)_{\alpha ,\beta }$ is
associative. If $G$ is commutative, then $(S^ME)_{\alpha ,\beta }$
is commutative. If $G$ is a Lie group, then $(S^ME)_{\alpha ,\beta
}$ is a Lie monoid. \par {\bf 3.} The $(S^ME)_{\alpha ,\beta }$ is
non-discrete, totally disconnected and infinite and non locally
compact for non degenerate $N$. Moreover, if $M$ and $N$ and $E$ are
dense in themselves, then the $(S^ME)_{\alpha ,\beta }$ is
topologically dense in itself and has the cardinality $card
[(S^ME)_{\alpha ,\beta }] \ge {\sf c} := card ({\bf Q_p})$.}

\par {\bf Proof.} We remind the following. Let $Q$ be a set and $T$
be a subset in $Q\times Q$. Then $T$ is called a relation on the set
$Q$. If $T$ satisfies Conditions $(E1-E3)$:
\par $(E1)$ $xTx$ for each $x\in Q$,
\par $(E2)$ from $xTy$ there follows $yTx$,
\par $(E3)$ $xTy$ and $yTz$ imply $xTz$,
then $T$ is called an equivalence relation.
\par We mention that each equivalence relation $T$ on $Q$ defines
some partition of $Q$ into non-intersecting subsets $A_x := \{ y\in
Q: xTy \} $ being equivalence classes in $Q$ relative to $T$. Thus
\par $(E4)$ $Q= \bigcup_{s\in S} A_s$ with $A_s\cap A_v = \emptyset $
for each $s\ne v\in S$, where $S$ is the corresponding set,
$S\subset Q$. Moreover, $x$ and $y\in A_s$ if and only if $xTy$.
Vice versa if a partition $ \{ A_s: s \in S \} $ of $Q$ into
pairwise disjoint subsets $A_s$ is given, $A_s\cap A_v=\emptyset $
for each $s\ne v\in S$, then it induces an equivalence relation $T$
on $Q$ so that $xTy$ if and only if $x$ and $y\in A_s$ (see also
\cite{eng}).
\par If $Q$ is a topological space and $T$ is some equivalence
relation on $Q$, then $Q/T$ denotes a set of all equivalence classes
in $Q$ relative to $T$. Then a mapping $q: Q\to Q/T$ exists posing
to each point $x\in Q$ its equivalence class $A_x$. This mapping $q$
is called the quotient mapping. In a class of all topologies on
$Q/T$ relative to which the quotient mapping is continuous a finest
exists: it is a family $\tau _{Q/T}$ of all subsets $U$ in $Q/T$ for
which $q^{-1}(U)$ is open. This topology $\tau _{Q/T}$ is called the
quotient topology. Moreover, $(Q/T, \tau _{Q/T})$ is called the
quotient space, $q$ is also called the natural quotient mapping or
shortly the natural mapping (see also \S 2.4 \cite{eng}).
\par Let $Q$ and $R$ be two topological spaces and let $f: Q\to R$
be  a continuous epimorphic mapping. It defines an equivalence
relation $T(f)$ on $Q$ generated by a partition $ \{ f^{-1}(y): y
\in R \} $. Then the mapping $f$ can be presented as a composition
$f = {\bar f}\circ q$, where $q: Q\to Q/T(f)$ is the natural
mapping, while $\bar f$ is a mapping from $Q/T(f)$ on $R$ prescribed
by the formula: ${\bar f}(f^{-1}(y)) =y$ for each $y\in R$.
Evidently $\bar f$ is the bijective continuous mapping, but
generally it need not be a homeomorphism. A continuous epimorphic
mapping $f: Q\to R$ is called a quotient mapping, if it is a
composition of a natural mapping $q: Q\to Q/T$ and some
homeomorphism, that is an equivalence relation $T$ on $Q$ exists and
a homeomorphism $h: Q/T\to R$ so that $f= h\circ q$. \par The
following proposition 2.4.3 \cite{eng} is useful. For a mapping $f$
of a topological space $A$ on a topological space $B$ the following
conditions are equivalent:
\par $(Q1)$ a function $f$ is a quotient mapping,
\par $(Q2)$ a set $f^{-1}(U)$ is  open in $A$ if and only if $U$ is
open in $B$,
\par $(Q3)$ a set $f^{-1}(C)$ is closed in $A$ if and only if $C$ is
closed in $B$,
\par $(Q4)$ a mapping ${\bar f} : A/T(f)\to B$ is a homeomorphism.
\par We remind that a set $\lambda $ is directed if it is supplied
with a relation $\le $ satisfying the following three conditions:
\par $(i)$ if $x\le y$ and $y\le z$, then $x\le z$;
\par $(ii)$ $x\le x$ for each $x\in \lambda $;
\par $(iii)$ for each pair $x, y \in \lambda $ an element $z\in \lambda $ exists
such that $x\le z$ and $y\le z$.
\par A subset $A$ in a directed set $\lambda $ is called cofinal if
for each $x\in \lambda $ an element $a\in A$ exists so that $x\le
a$. A set $\lambda $ is ordered if $(i,ii)$ are satisfied and the
following:
\par $(iv)$ if $x\le y$ and $y\le x$, then $x=y$.
\par An element $y$ of an ordered set $\lambda $ is called maximal
if from $y\le x$ the equality $x=y$ follows. If $\lambda $ is a set
and $\cal P$ is some property of its subsets, then $\cal P$ is of
finite character if the void set $\emptyset $ has it, and a subset
$A\subset \lambda $ possesses it if and only if each finite subset
of $A$ possesses it.
\par At first we consider trivial bundles with $G= \{ e \} $. So the
equivalence relation introduced in \S 2 we denote by $K_{\alpha
,\beta  }$ and it takes
the form: \\
$fK_{\alpha ,\beta }g$ if and only if there exist nets
$$\{ \psi _n\in Di^{\alpha }_{\beta ,0}(M) : \mbox{ }n\in \Omega \} ,$$
$$ \{ f_n\in C^{\alpha , w_0}_{\beta ,0} (M,N): \mbox{ } n \in \Omega  \}
\mbox{ and}$$
$$ \{ g_n\in C^{\alpha ,w_0}_{\beta ,0}(M,N): n \in \Omega  \}
\mbox{ such that}$$
$$(1)\mbox{ }f_n(x)=g_n(\psi _n(x))\mbox{ for each }x\in M\mbox{ and }
\lim_{n}f_n=f\mbox{ and }\lim_{n}g_n=g, $$ where $f, g \in C^{\alpha
,w_0}_{\beta ,0}( M,N)$ and the convergence is considered in this
space, $\Omega $ is a directed set.  Due to Condition $(1)$ these
equivalence classes $<f>_{K,\alpha ,\beta }$ are closed in
$C^{\alpha ,w_0}_{\beta ,0}( M,N)$. Then for $g\in <f>_{K,\alpha
,\beta }$ we write $gK_{\alpha , \beta }f$ also. The quotient space
$C^{\alpha ,w_0}_{\beta ,0} (M,N)/K_{\alpha ,\beta }$ we denote by
$(S^MN)_{\alpha ,\beta }$, where $w_0(M) = \{ y_0 \}$, $\theta
=w_0$.

\par Now we consider the wedge product $A\vee B := \rho ({\cal Z})$ be the wedge sum of
pointed spaces $(A,\{ a_{0,q}: q=1,...,k \} )$ and $(B,\{ b_{0,q}:
q=1,...,k \} )$, where ${\cal Z} := [A\times \{ b_{0,q}: q=1,...,k
\}\cup \{ a_{0,q}: q=1,...,k \}\times B]\subset A\times B$, $\rho $
is a continuous quotient mapping such that $\rho (x)=x$ for each
$x\in {\cal Z}\setminus \{ a_{0,q}\times b_{0,j}; q, j=1,...,k \} $
and $\rho (a_{0,q})=\rho (b_{0,q})$ for each $q=1,...,k$, where $A$
and $B$ are topological spaces with marked points $a_{0,q}\in A$ and
$b_{0,q}\in B$, $q=1,...,k$. Then the composition $g\circ f$ of two
elements $f, g\in C^{\alpha , w_0}_{\beta ,0}(M, N)$ is defined on
the domain \par $(W1)$ ${\bar M}\vee {\bar M}\setminus \{
s_{0,q}\times s_{0,q}: q=1,...,k \} =: M\vee M$.
\par Let $M= {\bar M}\setminus \{ s_{0,q}: q=1,..., k \} $ be as in \S 1.
In view of Conditions 2.4$(D1-D8)$ we can choose a refinement of an
initial covering.
\par We shall use the Teichm\"uller-Tukey's lemma \cite{eng}.
If $\lambda $ is a set, while $\cal P$ is a property of a finite
order, then each subset $A\subset \lambda $ having this property
$\cal P$ is contained in a set $B$ also satisfying $\cal P$ and $B$
is a maximal element in an ordered by inclusion family of all
subsets of $\lambda $ having the property $\cal P$.
\par  The topological space $(M, \tau _{p,M})$ is totally disconnected and it is not compact.
By its construction the set $\breve{M}$ is open in $\bar M$, since
$\breve{M} = \bigcup \{ V: ~ h_{\bar M}(V)$ $\mbox{is open in}$ $X$
$\mbox{for some plot}$ $h\in {\cal P}_{\bar M} \} $, where ${\cal
P}_{\bar M}$ is a family of plots defining a $C^{\alpha '}_{\beta }$
differentiable structure of $\bar M$ (see also \S 3.1 above). Put
$A_x = \{ x \} $ for each $x\in \breve{M}\setminus \{ s_{0,q}:
q=1,...,k \} $ and $\bigcup_{q=1}^k A_{s_{0,q}} = {\bar M}\setminus
(\breve{M}\setminus \{ s_{0,q}: q=1,...,k \} )$, each $A_{s_{0,q}}$
is clopen in ${\bar M}\setminus (\breve{M}\setminus \{ s_{0,q}:
q=1,...,k \} )$, $A_{s_{0,q}}\cap A_{s_{0,t}}=\emptyset $ for each
$1\le q\ne t\le k$. Such partition $ \{ A_x \} $ of $\bar M$ induces
an equivalence relation $T$ in $\bar M$ (see above). Using the
family of (generalized) diffeomorphisms $Di^{\alpha }_{\beta ,0}(M)$
and the quotient space ${\bar M}/T$ in case of necessity we can
reduce our proof to the case, when classes of equivalent mappings
are considered on $\bar M$ satisfying the condition \par $(W2)$
${\bar M}\setminus (\breve{M}\setminus \{ s_{0,q}: q=1,...,k \} ) =
\{ s_{0,q}: q=1,...,k \} .$
\par We consider families which are bases of $\tau _{p,{\bar M}}$ open
neighborhoods $W_{v,q}$ of marked points $s_{0,q}$ in $\bar M$, that
is $s_{0,q}\in W_{v,q}$. Each family $ \{ W_{v,q}: v\in \lambda
_{q,M} \} $ is ordered by inclusion: $W_{u,q}\le W_{v,q}$ if and
only if $W_{v,q}\subset W_{u,q}$. Each finite intersection of open
sets is open. So we fix an infinite (generalized) atlas
\par $(2)$ ${\tilde A}t'(M):= \{ ({{\tilde U}'}_j,{\phi '}_j): j\in \lambda \}$
such that ${\phi '}_j: {{\tilde U}'}_j\to B_j$ are homeomorphisms on
bounded ${\cal A}_r$ convex subsets $B_j$ in $X$, each ${{\tilde
U}'}_j$ is clopen in $(M,\tau _{p,M})$, $\lambda $ is a directed
set. Moreover,
\par $(3)$ for each $q\in \{ 1,..., k \} $ and every $u\in \lambda
_{q,M}$ there exists $v\in \lambda $ so that for each $v\le j \in
\lambda $ either
\par ${{\tilde U}'}_j\subset W_{u,q}$ or ${{\tilde U}'}_j\cap W_{u,q} =
\emptyset $. Then also \par $(4)$ for each $q$ a subset $\omega
_q\subset \lambda $ exists so that $cl_{\bar M}[\bigcup_{j\in \omega
_q}{{\tilde U}'}_j]$ is a clopen neighborhood of $s_{0,q}$ in $\bar
M$, where $cl_{\bar M}A$ denotes the closure of a subset $A$ in
$({\bar M}, \tau _{p,{\bar M}})$. Property $(4)$ follows from
$(2,3)$. The topological space $(M,\tau _{p,M})$ is not compact,
hence
\par $(5)$ $card (\lambda )\ge card ({\bf Z})=\aleph _0$.
\par  In the wedge product
$M\vee M$ we choose the following atlas \par $(6)$ ${\tilde
A}t'(M\vee M)=\{ (W_l, \xi _l): l\in \mu \}$ such that $\xi _l:
W_l\to C_l$ are homeomorphisms, $C_l = C_{l,1}\vee ' C_{l,2}$, each
$C_{l,1}$ and $C_{l,2}$ are bounded ${\cal A}_r$ convex subsets in
$X$, where we denote $[(C_{l,1}\cup \{ x_{0,q}: q=1,...,k \})\vee
(C_{l,2}\cup \{ x_{0,q}: q=1,...,k \})]\setminus \{ x_{0,q}\times
x_{0,q}: q=1,...,k \}  =: C_{l,1}\vee ' C_{l,2}$ for suitable
distinct marked points $x_{0,q}$ in $X$ corresponding to $s_{0,q}$,
$\mu $ is a directed set; also
\par $(7)$ for each $q_1, q_2 \in \{ 1,..., k \} $ and every $u_1\in \lambda _{q_1,M}$
and each $u_2 \in \lambda _{q_2,M}$ there exists $v\in \mu $ so that
for each $v\le l \in \mu $ either
\par $W_l\subset W_{u_1,q_1}\vee W_{u_2,q_2}$ or $W_l\cap (W_{u_1,q_1}\vee W_{u_2,q_2}) =
\emptyset $. Therefore, \par $(8)$ for each $q_1, q_2\in \{ 1,...,k
\} $ a subset $\nu _{q_1,q_2}\subset \mu $ exists so that $cl_{{\bar
M}\vee {\bar M}} [\bigcup_{l\in \nu _{q_1,q_2}} W_l]$ is a clopen
neighborhood of $s_{0,q_1}\vee s_{0,q_2}$ in ${\bar M}\vee {\bar
M}$, where $cl_{{\bar M}\vee {\bar M}}A$ denotes the closure of a
subset $A$ in $({\bar M}\vee {\bar M}, \tau _{p,{\bar M}}\times \tau
_{p,{\bar M}})$. We get property $(8)$ from Conditions $(6,7)$.
Since the topological space $(M\vee M,\tau _{p,M}\times \tau
_{p,M})$ also is not compact (see $(W1,W2)$ above), then the
cardinalities of $\lambda $ and $\mu $ are the same, so we can
choose $\omega _q$ and $\nu _{q_1,q_2}$ so that
\par $(9)$ $card (\lambda \setminus  \bigcup_{q=1}^k
\omega _q ) = card (\lambda \setminus  \omega _{q_1} ) = card (\mu
\setminus \bigcup_{q_1,q_2\in \{ 1,...,k \} } \nu _{q_1,q_2} ) =
card (\mu \setminus \nu _{q_1,q_2}) \ge \aleph _0$, also $card (\nu
_{q_1,q_2}) = card (\omega _{q_1})= card (\omega _{q_2}) \ge \aleph
_0$ for all $q_1, q_2\in \{ 1,...,k \} $, where the product topology
$\tau _{p,M}\times \tau _{p,M}$ on $M\times M$ induces the
corresponding topology on the subset $M\vee M$. We denote this
topology on $M\vee M$ by $\tau _{p,M}\times \tau _{p,M}$ or $\tau
_{p,M\vee M}$. Due to the Zermelo's Theorem \cite{eng} we can
consider sets $\omega _q$ and $\nu _{q_1,q_2}$ as ordinals $[\omega
_q]$ and $[\nu _{q_1,q_2}]$ of the same type $[\omega _{q_1}] =
[\omega _{q_2}] = [\nu _{q_1,q_2}]$ for all $q_1, q_2\in \{ 1,...,k
\} $.
\par In view of the Teichm\"uller-Tukey's lemma we can choose
$W_l$ consistent with ${{\tilde U}'}_{j_1}\vee  {{\tilde U}'}_{j_2}$
such that to fix a (generalized) $C^{\alpha '}_{\beta
}$-diffeomorphisms $\chi : M\vee M\to M$ satisfying the following
conditions $(10-12)$:
$$(10)\mbox{ }\chi (W_l)=
{{\tilde U}'}_l\mbox{ for each }l\in [\nu _{q_1,q_2}]~ \forall q_1,
q_2 \in \{ 1,...,k \} \mbox{ and}$$
$$(11)\mbox{ }\chi (W_l)=
{{\tilde U}'}_{\kappa (l)} \mbox{ for each }l\in (\mu \setminus
\bigcup_{q_1, q_2\in \{ 1,...,k \} } \nu _{q_1,q_2}),\mbox{ where}$$
$$(12)\mbox{ }\kappa : (\mu \setminus \bigcup_{q_1, q_2\in \{ 1,...,k \} } \nu
_{q_1,q_2}) \to (\lambda \setminus \bigcup_{q\in \{ 1,...,k \} }
\omega _q)$$ is a bijective mapping (see also \S 2 above). This
induces the continuous injective homomorphism
$$(13)\mbox{ }\chi ^*: C^{\alpha ,w_0}_{\beta ,0}(
(M\vee M, \{ s_{0,q_1}\times s_{0,q_2}: q_1, q_2 \in \{ 1,...,k \}
);(N,y_0))$$  $$\to C^{\alpha ,w_0}_{\beta ,0} ((M, \{ s_{0,q}:
q=1,...,k \} ); (N,y_0))\mbox{ such that }$$
$$(14)\mbox{ }\chi ^*(g\vee f)(x)=(g\vee f)(\chi ^{-1}(x))$$
for each $x\in M$, where $(g\vee f) (y)=f(y)$ for each $y\in M_2$
and $(g\vee f)(y)=g(y)$ for every $y\in M_1$, $M_1\vee M_2=M\vee M$,
$M_i=M$ for $i=1,2$ are two copies of $M$. Therefore
$$(15)\mbox{ }g\circ f := \chi ^*(g\vee f)$$
may be considered as defined on $M$ also, that is, to $g\circ f$
there corresponds the unique element in $C^{\alpha ,w_0}_{\beta ,0}
((M, \{ s_{0,q}: q=1,...,k \} ); (N,y_0))$. \par We have $f(\psi
)\in C^{\alpha ,w_0}_{\beta ,0}((M,\{ s_{0,q}: q=1,...,k \} );
(N,y_0))$ for each $f\in C^{\alpha , w_0}_{\beta ,0}((M,\{ s_{0,q}:
q=1,...,k \}  ); (N,y_0))$ and $\psi \in Di^{\alpha }_{\beta ,0}(M)$
due to Lemma 9 and Corollary 10 in \cite{ludfnafjms09} applied
uniformly by finite dimensional over $\bf K$ embedded into $\bar M$
differentiable subspaces with the corresponding (generalized)
atlases (see also \cite{lutmf99,luum985}). The diffeomorphism $\chi
: M\vee M\to M$ is of class $C^{\alpha '}_{\beta }$, $\alpha '\ge
\alpha $, and from Conditions 2$(4,5)$ for $f_i\in C^{\alpha
,w_0}_{\beta ,0}((M, \{ s_{0,q}: q=1,...,k \} ); (N,y_0))$ it
follows that for $f=\chi ^*(f_1\vee f_2)$ also Condition 2$(4)$ is
satisfied, since $\chi $ fulfils Conditions $(10-14)$. Moreover, ~
$<f>_{K,\alpha ,\beta }\circ <g>_{K,\alpha ,\beta }= <f\vee
g>_{K,\alpha ,\beta }$ for each $f$ and $g \in C^{\alpha
,w_0}_{\beta ,0} ((M, \{ s_{0,q}: q=1,...,k \} ); (N,y_0))$, since
if ${\tilde f}_n(x)= f_n(\eta _n(x))$ and ${\tilde g}_n(x)=g_n(\zeta
_n(x))$ for each $x\in M$, then $({\tilde f}_n\vee {\tilde
g}_n)(x)=(f_n(\eta _n)\vee g_n( \zeta _n))(x),$ where $\eta _n$ and
$\zeta _n \in Di^{\alpha }_{\beta ,0} (M)$. Hence the composition is
continuous for the quotient space. \par In view of Conditions
2$(4,5)$ for each $f\in C^{\alpha ,w_0}_{\beta ,0} ((M, \{ s_{0,q}:
q=1,...,k \} ); (N,y_0))$ there exist nets $\{ \psi _n: n \in \Omega
\}$, $\{ \eta _n: n \in \Omega \}$ and $\{ \zeta _n: n \in \Omega \}
$ in $Di^{\alpha }_{\beta ,0} (M)$, $\{ f_n: n \in \Omega \} $, $\{
w_{0,n}: n \in \Omega \} $ and $\{ g_n: n \in \Omega \} $ in
$C^{\alpha ,w_0}_{\beta ,0} ((M, \{ s_{0,q}: q=1,..., k \} );
(N,y_0))$ such that $[w_{0,n}\vee f_n] (\psi _n\vee \eta _n (x)) =
g_n(\zeta _n(\chi (x)))$, where
$$(17)\mbox{ }\lim_{n}f_n=f,\mbox{ }\lim_{n}g_n=g,$$
$$(18)\mbox{ }\lim_{n}w_{0,n}=w_0,\mbox{ }
f_n(x)\ne y_0\mbox{ for each }x\in M,$$  $\Omega $ is a directed
set. On the other hand, from $\lim_{n}(f_n\vee g_n)=f\vee g$ it
follows that $\lim_{n}f_n=f$ and $\lim_{n}g_n=g$. We choose
\par $(19)$ $ \{ \zeta _n: n \in \Omega \} $ so that for each open
neighborhood $U$ of $ \{ s_{0,q}: q=1,...,k \} $ in $({\bar M}, \tau
_{p,{\bar M}})$ there exists $m\in \Omega $ for which $\zeta _n(\chi
([M\times \{ s_{0,q} : q=1,...,k \} ] \cap [M\vee M])) \subset U$
for every $n\ge m$ in $\Omega $. \par Using Formulas $(10-15)$ and
Conditions 2$(4-6)$ we get $<w_0\circ f>_{K,\alpha ,\beta
}=<f>_{K,\alpha ,\beta }$ and $<w_0>_{K,\alpha ,\beta }=e$ is the
unit element in $(S^MN)_{\alpha ,\beta }$, since $<f>_{K,\alpha
,\beta }\circ <g>_{K,\alpha ,\beta }=<f\vee g>_{K,\alpha ,\beta }$
for each $f$ and $g \in C^{\alpha ,w_0}_{\beta ,0} ((M, \{ s_{0,q}:
q=1,...,k \} ); (N,y_0))$.

\par Let us consider now the general case of fiber bundles.
If a homomorphism $\theta : G\to F$ of $C^{\alpha '}_{\beta
}$-differentiable groups exists, then an induced principal $F$ fiber
bundle $(E\times ^{\theta }F)(N,F,\pi ^{\theta }, \Psi ^{\theta })$
is given with the total space $(E\times ^{\theta }F)= (E\times
F)/{\cal Y}$, where $\cal Y$ is the equivalence relation such that
$(vg,f){\cal Y} (v,\theta (g)f)$ for each $v\in E$, $g\in G$, $f\in
F$. Then the projection $\pi ^{\theta }: (E\times ^{\theta }F)\to N$
is defined by the equation $\pi ^{\theta }([v,f]) =\pi (v)$, where
$[v,f] := \{ (w,b): (w,b){\cal Y} (v,f), w\in E, b\in F \} $ denotes
the equivalence class of $(v,f)$. \par This implies that each
parallel transport structure $\bf P$ on the principal $G$ fiber
bundle $E(N,G,\pi ,\Psi )$ induces a parallel transport structure
${\bf P}^{\theta }$ on the induced bundle by the formula ${\bf
P}^{\theta }_{{\hat \gamma }, [u,f]}(x) = [{\bf P}_{{\hat \gamma
},u}(x),f]$.

\par We define multiplication with the help of certain embeddings and
isomorphisms of spaces of functions. Mention that for each two
$C^{\alpha '}_{\beta }$ diffeomorphic $C^{\alpha '}_{\beta }$
differentiable spaces $A$ and $B$ in a topological vector space $X$
over ${\cal A}_r^l$ the spaces $C^{\alpha }_{\beta }(A,Y)$ and
$C^{\alpha }_{\beta }(B,Y)$ are isomorphic as topological ${\cal
A}_r$ vector spaces, where $Y$ is also a topological vector space
over ${\cal A}_r$, $\alpha \le \alpha '$, consequently, $C^{\alpha
}_{\beta }(A,N)$ and $C^{\alpha }_{\beta }(B,N)$ are isomorphic as
uniform spaces. Naturally we consider the space \par $C^{\alpha
}_{\beta }(M, \{ s_{0,1},...,s_{0,k} \} ; {\cal W},y_0) := \{ (E,f):
E=E(N,G,\pi ,\Psi )\in {\cal W}, f={\bf P}_{{\hat \gamma },y_0}\in
C^{\alpha }_{\beta }: \pi \circ f(s_{0,q})=y_0 \forall q=1,...,k;
\pi \circ f={\hat \gamma }, \gamma \in C^{\alpha ,w_0}_{\beta
,0}(M,N) \} $ \\ which is the space of all $C^{\alpha '}_{\beta }$
principal $G$ fiber bundles $E$ with their parallel transport
$C^{\alpha }_{\beta }$-mappings $f={\bf P}_{{\hat \gamma },y_0}$ in
accordance with \S 2, where $\cal W$ is the same family of all
principal $C^{\alpha '}_{\beta }$-fiber bundles $E(N,G,\pi ,\Psi )$
as in \S 2.6. Put $\omega _0=(E_0,{\bf P}_0)$ be its element such
that $\gamma _0(M)= \{ y_0 \} $, where $e\in G$ denotes the unit
element, $E_0=N\times G$, $\pi _0(y,g)=y$ for each $y\in N$, $g\in
G$, ${\bf P}_{{\hat \gamma }_0,u}={\bf P}_0$.
\par The mapping $\Xi : {\hat M}\to M$ from \S 1 induces the embedding
\par $\Xi ^*: C^{\alpha }_{\beta }(M,\{ s_{0,1},...,s_{0,k} \} ;{\cal W},y_0)\hookrightarrow
C^{\alpha }_{\beta }({\hat M},\{ {\hat s}_{0,1},...,{\hat s}_{0,2k}
\} ;{\cal W},y_0)$. We consider the wedge product $g\vee f$ of two
elements $f, g\in C^{\alpha }_{\beta }((M, \{ s_{0,1},...,s_{0,k}
\}) ;(N,y_0))$ which is defined on the domain $M\vee M$, where to
$f, g$ two mappings $f_1, g_1 \in C^{\alpha }_{\beta }(({\hat M}, \{
{\hat s}_{0,1},...,{\hat s}_{0,2k} \} ); (N,y_0))$ correspond such
that $f_1 = f\circ \Xi $ and $g_1 = g\circ \Xi $.

\par Suppose that $(E_j,{\bf P}_{{\hat \gamma }_j,u^j})\in
C^{\alpha }_{\beta }(M,\{ s_{0,1},...,s_{0,k} \} ;{\cal W},y_0)$,
$j=1, 2$, then we take their wedge product ${\bf P}_{{\hat \gamma
},u^1} := {\bf P}_{{\hat \gamma }_1,u^1}\vee {\bf P}_{{\hat \gamma
}_2,v}$ on $M\vee M$ with $v_q=u_q g_{2,q}^{-1}g_{1,q+k}= y_0\times
g_{1,q+k}$ for each $q=1,...,k$ due to the alternativity of $G$,
$\gamma =\gamma _1\vee \gamma _2$, where ${\bf P}_{{\hat \gamma
}_j,u^j}({\hat s}_{j,0,q})= y_0\times g_{j,q}\in E_{y_0}$ for every
$j$ and $q$. For each $\gamma _j: M\to N$ a mapping ${\tilde \gamma
}_j: M\to E_j$ exists such that $\pi \circ {\tilde \gamma }_j=\gamma
_j$. We denote by ${\bf m}: G\times G\to G$ the multiplication
operation in the group $G$. The wedge product $(E_1,{\bf P}_{{\hat
\gamma }_1,u^1})\vee (E_2,{\bf P}_{{\hat \gamma }_2,u^2})$ is the
principal $G$ fiber bundle $(E_1\times E_2)\times ^{\bf m}G$ with
the parallel transport structure ${\bf P}_{{\hat \gamma }_1,u^1}\vee
{\bf P}_{{\hat \gamma }_2,v}$.

\par Due to Conditions $(10-14)$ and 2$(4,5)$ we get
the following embedding $\chi ^*: C^{\alpha }_{\beta }(M\vee M,\{
s_{0,q}\times s_{0,q}: q=1,...,k \} ;{\cal W},y_0)\hookrightarrow
C^{\alpha }_{\beta }(M,\{ s_{0,q}: q=1,...,k \} ;{\cal W},y_0)$.
Therefore, $g\circ f := \chi ^*(f\vee g)$ is the composition in
$C^{\alpha }_{\beta }(M,\{ s_{0,q}: q=1,...,k \} ;{\cal W},y_0)$.
\par Generalizing the beginning of this section we define
the following equivalence relation $K_{\alpha ,\beta }$ in
$C^{\alpha }_{\beta }(M,\{ s_{0,q}: q=1,...,k \} ;{\cal W},y_0)$:
$fK_{\alpha ,\beta }h$ if and only if nets $\eta _n\in Di^{\alpha
}_{\beta ,0}(M)$, also $f_n$ and $h_n\in C^{\alpha }_{\beta }(M,\{
s_{0,q}: q=1,...,k \} ;{\cal W},y_0)$ with $\lim_n f_n=f$ and
$\lim_n h_n=h$ such that $f_n(x)=h_n(\eta _n(x))$ for each $x\in M$
and $n\in \omega $, where $\omega $ is a directed set and
convergence is considered in $C^{\alpha }_{\beta }(M,\{ s_{0,q}:
q=1,...,k \} ;{\cal W},y_0)$.
\par Thus the following quotient uniform space \\
$C^{\alpha }_{\beta }(M,\{ s_{0,q}: q=1,...,k \} ; {\cal
W},y_0)/K_{\alpha ,\beta } =: (S^ME)_{\alpha ,\beta }$ exists.

\par We consider an element $f={\bf P}_{{\hat \gamma
},u}$ as $f\circ \Xi ^{-1}$ on ${\bar M}\setminus {\bar M}_f$, where
$\pi \circ f = {\hat \gamma }$, ${\hat \gamma }=\gamma \circ \Xi $.
We denote $f\circ \Xi ^{-1}$ also by $f$. If $M$ and $N$ and $G$ are
separable, then $C^{\alpha }_{\beta }(M,\{ s_{0,q}: q=1,...,k \}
;{\cal W},y_0)$ is separable, consequently, $(S^ME)_{\alpha ,\beta
}$ is also separable.

\par By our construction each equivalence class $<f>_{K,\alpha
,\beta  }$ is closed in $C^{\alpha }_{\beta }(M, \{ s_{0,q}:
q=1,...,k \} ; {\cal W},y_0)$. Therefore, each point $g$ in
$(S^ME)_{\alpha ,\beta }$ is closed in it. A topological space $S$
is $T_1$ if and only if each singleton (one-pointed set) $ \{ g \} $
is closed in it (see \S 1.5 \cite{eng}). Thus the topological space
$(S^ME)_{\alpha ,\beta }$ possesses the $T_1$ separability axiom.

\par The uniform space $C^{\alpha }_{\beta }(M,\{ s_{0,q}: q=1,...,k \} ;{\cal W},y_0)$
is complete due to Theorem 12.1.4 \cite{nari}, when $N$ and $G$ are
complete. Each class of $K_{\alpha ,\beta }$-equivalent elements is
closed in it. Consider reparametrizations of elements $f$ of
$C^{\alpha }_{\beta }(M,\{ s_{0,q}: q=1,...,k \} ;{\cal W},y_0)$
relative to the action $f\mapsto f\circ \psi $, $\psi \in Di^{\alpha
}_{\beta ,0}(M)$, of the family $Di^{\alpha }_{\beta ,0}(M)$ on
$\bar M$. Then to each Cauchy net in $(S^ME)_{\alpha ,\beta }$ there
corresponds a Cauchy net in $C^{\alpha }_{\beta }(M, \{ s_{0,q}:
q=1,...,k \} ; {\cal W},y_0)$. Hence $(S^ME)_{\alpha ,\beta }$ is
complete, if $N$ and $G$ are complete.

\par If $f, g\in C^{\alpha }_{\beta } (M,X)$ and $f(M)\ne g(M)$, then
\par $< f\circ
\psi - g >_{K,\alpha ,\beta } \ne <w_0\times e >_{K,\alpha ,\beta }$
for each $\psi \in Di^{\alpha }_{\beta ,0}(M)$. Thus equivalence
classes $<f>_{K,\alpha ,\beta }$ and $<g>_{K,\alpha ,\beta }$ are
different.
\par The uniform space $C^{\alpha }_{\beta }(M, \{
s_{0,q}: q=1,...,k \} ; {\cal W},y_0)$ is totally disconnected and
dense in itself, when $M$ and $N$ and $E$ are dense in themselves,
since $C^{\alpha }_{\beta }(M,Y)$ is such for each topological
vector space $Y$ over ${\cal A}_r$. Thus, the uniform space
$(S^ME)_{\alpha , \beta }$ is non-discrete and dense in itself.
\par Take a restriction of $E$ for a chosen (generalized) chart $(\phi ^U, U)$ of $E$.
In accordance with conditions on $E$ there exists $u\in U$ such that
$V := \phi ^U(U) - u$ is absolutely ${\cal A}_r$ convex, where $X_E$
is a topological ${\cal A}_r$ vector space on which $E$ is modeled,
$\phi ^U: U\to X_E$. The subspace $X_{E,u} := cl_{X_E} (span_{{\cal
A}_r} V)$ we call the (generalized) tangent space at $u\in U$ to
$E$. Then we can choose $U$ so that $\pi (U) =: U_N$ and $(\phi
^U_N, U_N)$ is the chart of $N$, hence $X_{N,y} := cl_{X_N}
(span_{{\cal A}_r} V_N)$ is the (generalized) tangent space at $y =
\pi (u) \in U_N$ to $N$, where $V_N := \phi ^U_N(U_N) - y$ is
absolutely ${\cal A}_r$ convex in a topological vector space $X_N$
on which $N$ is modeled, $\pi : E\to N$ is the projection mapping of
the fiber bundle. For non degenerate $N$ the space $X_{N,y}$ has a
dimension over ${\cal A}_r$ not less than one. When $\alpha '\ge
\alpha +1$ the tangent bundle $TC^{\alpha }_{\beta }(M,E_U)$ is
isomorphic with $C^{\alpha }_{\beta }(M,TE_U)$, where $TE_U$ is the
$C^{\alpha '-1}_{\beta }$ fiber bundle. There is an infinite family
of $f_{\alpha }\in C^{\alpha }_{\beta }(M,TE_U)$ with pairwise
distinct images in $TE_U$ for different $\alpha $ such that
$f_{\alpha }(M)$ is not contained in $\bigcup_{\beta <\alpha
}f_{\beta }(M)$, $\alpha \in \Lambda $, where $\Lambda $ is an
infinite ordinal. Therefore, $T(S^ME_U)_{\alpha ,\beta }$ is an
infinite dimensional fiber bundle due to $(ii)$. We say that
$(S^ME)_{\alpha ,\beta }$ is infinite dimensional over ${\cal A}_r$
if for each (generalized) chart $U$ of $E$ the fiber bundle
$T(S^ME_U)_{\alpha ,\beta }$ is infinite dimensional over ${\cal
A}_r$. For two fiber bundles $E_1$ and $E_2$ which are $C^{\alpha
}_{\beta }$ isomorphic the uniform spaces $(S^ME_1)_{\alpha ,\beta
}$ and $(S^ME_2)_{\alpha ,\beta }$ are $C^{\alpha }_{\beta }$
isomorphic. Thus $(S^ME)_{\alpha ,\beta }$ is infinite and non
locally compact for each $\alpha '\ge \alpha $, since
$(S^ME_U)_{\alpha ,\beta }$ is infinite dimensional over ${\cal
A}_r$ and there is an embedding $(S^ME_U)_{\alpha ,\beta
}\hookrightarrow (S^ME)_{\alpha ,\beta }$.
\par Evidently, if $f\vee g=h\vee g$ or $g\vee f=g\vee h$ for
$\{ f, g, h \} \subset C^{\alpha }_{\beta }(M,\{ s_{0,q}: q=1,...,k
\} ;{\cal W},y_0)$, then $f=h$. Thus $\chi ^*(f\vee g)=\chi ^*(h\vee
g)$ or $\chi ^*(g\vee f)=\chi ^*(g\vee h)$ is equivalent to $f=h$
due to the definition of $f\vee g$ and the definition of equal
functions, since $\chi ^*$ is the embedding. Using the equivalence
relation $K_{\alpha ,\beta }$ gives $<f>_{K,\alpha ,\beta }\circ
<g>_{K,\alpha ,\beta }= <h>_{K,\alpha ,\beta }\circ <g>_{K,\alpha
,\beta }$ or $<g>_{K,\alpha ,\beta }\circ <f>_{K,\alpha ,\beta }=
<g>_{K,\alpha ,\beta }\circ <h>_{K,\alpha ,\beta }$ is equivalent to
$<h>_{K,\alpha ,\beta }=<f>_{K,\alpha ,\beta }$. Therefore,
$(S^ME)_{\alpha ,\beta }$ has the cancelation property.
\par The group $G$ is alternative, so
$a_{2,q}[a_{2,q}^{-1}(a_{2,q+k}(a_{2,q}^{-1}a_{1,q+k}))]=
a_{2,q+k}(a_{2,q}^{-1}a_{1,q+k})$, hence ${\bf P}_1\vee ({\bf
P}_2\vee {\bf P}_2)= ({\bf P}_1\vee {\bf P}_2)\vee {\bf P}_2$; also
$a_{2,q}[a_{2,q}^{-1}(a_{1,q+k}(a_{1,q}^{-1}a_{1,q+k}))]=
a_{1,q+k}(a_{1,q}^{-1}a_{1,q+k})$, consequently, ${\bf P}_1\vee
({\bf P}_1\vee {\bf P}_2)= ({\bf P}_1\vee {\bf P}_1)\vee {\bf P}_2$
and inevitably for equivalence classes $(aa)b=a(ab)$ and
$b(aa)=(ba)a$ for each $a, b\in (S^ME)_{\alpha ,\beta }$. Thus the
$(S^ME)_{\alpha ,\beta }$ is alternative. \par Evidently $M\vee
(M\vee M)$ is (generalized) $C^{\alpha }_{\beta ,0}$-diffeomorphic
with $(M\vee M)\vee M$ (see 2$(6)$).
\par If $G$ is associative, then the parallel transport structure
gives $(f\vee g)\vee h=f\vee (g\vee h)$ on $M\vee M\vee M$ for each
$\{ f, g, h \} \subset C^{\alpha }_{\beta }(M,\{ s_{0,q}: q=1,...,k
\} ;{\cal W},y_0)$. Applying the embedding $\chi ^*$ and the
equivalence relation $K_{\alpha ,\beta }$ we get, that
$(S^ME)_{\alpha ,\beta }$ is associative $<f>_{K,\alpha ,\beta
}\circ (<g>_{K,\alpha ,\beta }\circ <h>_{K,\alpha ,\beta })=
(<f>_{K,\alpha ,\beta  }\circ <g>_{K, \alpha , \beta })\circ <h>_{K,
\alpha , \beta }$.

\par Let $w_0$ be a mapping $w_0: M\to W$ such that $w_0(M)=\{ y_0
\times e \} $. Consider $w_0\vee (E,f)$ for some $(E,f)\in C^{\alpha
}_{\beta } (M,\{ s_{0,q}: q=1,...,k \} ;{\cal W},y_0)$. Suppose
$(E,f)\in C^{\alpha }_{\beta }(M,\{ s_{0,q}: q=1,..,k \}; {\cal
W},y_0)$. A net $U_n$ of open or canonical closed subsets in $M$
exists such that $\bigcap_n U_n = \{ s_{0,q}: q=1,...,k \} $. We
mention that for each marked point $s_{0,q}$ in $M$ there exists a
neighborhood $U$ of $s_{0,q}$ in $M$ such that for each $\gamma
_1\in C^{\alpha }_{\beta }((M, \{ s_{0,q}: q=1,...,k \} ); (N,y_0))$
there exists $\gamma _2\in C^{\alpha }_{\beta }$ such that they are
$K_{\alpha ,\beta }$ equivalent and $\gamma _2|_U = y_0$ due to
Conditions 2$(4-6)$ and $(10-15)$.
\par A net $\eta _n \in Di^{\alpha }_{\beta ,0}(M)$ exists
together with $w_n, f_n\in C^{\alpha }_{\beta }(M,\{ s_{0,q}:
q=1,...,k \} ;{\cal W},y_0)$ with
\par  $(20)$ $\lim_{n}f_n=f$, $\lim_{n
}w_n=w_0$ and $\lim_{n}\chi ^*(f_n\vee w_n)(\eta _n^{-1})=f$ due to
$\pi \circ f(s_{0,q})=s_{0,q}$ in the formula of difference
quotients of compositions of functions (see also $(17-19)$ above).
Indeed, we can apply Lemma 9 and Corollary 10 in \cite{ludfnafjms09}
uniformly by finite dimensional over $\bf K$ embedded into $\bar M$
differentiable subspaces with the corresponding (generalized)
atlases.

\par Therefore, $w_0\vee (E,f)$ and
$(E,f)$ belong to the equivalence class $<(E,f)>_{K,\alpha ,\beta }
:= \{ g\in C^{\alpha }_{\beta }(M,\{ s_{0,q}: q=1,...,k \} ;{\cal
W},y_0): (E,f)K_{\alpha ,\beta }g \} $ due to $(20)$. Thus,
$<w_0>_{K,\alpha ,\beta }\circ <g>_{K,\alpha ,\beta }= <g>_{K,\alpha
,\beta }$.
\par  The $C^{\alpha '}_{\beta }$ differentiable space $(M\vee
M)\setminus \{ s_{0,q}\times s_{0,j}: q, j=1,...,k \} $ is open in
${\bar M}\vee {\bar M}$ and has the $C^{\alpha }_{\beta
}$-diffeomorphism $\psi $ such that $\psi (x,y)=(y,x)$ for each
$(x,y)\in ((M\times M)\setminus \{ s_{0,q}\times s_{0,j}: q,
j=1,...,k \} )$. Suppose now, that $G$ is commutative. Then $(f\vee
g)\circ \psi |_{(M\times M\setminus \{ s_{0,q}\times s_{0,j}: q,
j=1,...,k \} )}= g\vee f|_{(M\times M\setminus \{ s_{0,q}\times
s_{0,j}: q, j=1,...,k \} )}$. On the other hand, $<f\vee
w_0>_{K,\alpha ,\beta }=<f>_{K,\alpha ,\beta }=<f>_{K,\alpha ,\beta
}\circ <w_0>_{K,\alpha ,\beta }=<w_0>_{K,\alpha ,\beta }\circ
<f>_{K,\alpha ,\beta }$, hence, $<f\vee g>_{K,\alpha ,\beta
}=<f>_{K,\alpha ,\beta }\circ <g>_{K,\alpha ,\beta }=<f\vee
w_0>_{K,\alpha ,\beta }\circ <w_0\vee g>_{K,\alpha ,\beta  }=<(f\vee
w_0)\vee (w_0\vee g)>_{K,\alpha ,\beta }= <(w_0\vee g)\vee (f\vee
w_0)>_{K,\alpha ,\beta }$ due to the existence of the unit element
$<w_0>_{K,\alpha ,\beta }$ and due to the properties of $\psi $.
Indeed, take a net $\psi _n$ as above. Therefore, the parallel
transport structure gives $(g\vee f)(\psi (x,y))=(g\circ f)(y,x)$
for each $x, y\in M$, consequently, $(f\circ g)K_{\alpha ,\beta
}(g\circ f)$ for each $f, g\in C^{\alpha }_{\beta }(M, \{ s_{0,q}:
q=1,...,k \} ;{\cal W},y_0)$. The using of the embedding $\chi ^*$
gives that $(S^ME)_{K,\alpha ,\beta }$ is commutative, when $G$ is
commutative.
\par The mapping $(f,g)\mapsto f\vee g$ from $C^{\alpha }_{\beta ,0}(M,\{ s_{0,q}:
q=1,...,k \} ;{\cal W},y_0)^2$  into $C^{\alpha }_{\beta ,0}(M\vee
M\setminus \{ s_{0,q}\times s_{0,j}: q, j=1,...,k \}; {\cal W},y_0)$
is of class $C^{\alpha }_{\beta }$. Since the mapping $\chi ^*$ is
of class $C^{\alpha }_{\beta }$, then $(f,g)\mapsto \chi ^*(f\vee
g)$ is the $C^{\alpha }_{\beta }$-mapping. The quotient mapping from
$C^{\alpha }_{\beta }(M,\{ s_{0,q}: q=1,...,k \} ;{\cal W},y_0)$
into $(S^ME)_{K,\alpha ,\beta }$ is continuous and induces the
quotient uniformity. On the other hand, $T^b(S^ME)_{K,\alpha ,\beta
}$ has embedding into $(S^MT^bE)_{\alpha ,\beta }$ for each $1\le
b\le \alpha '-\alpha $, when $\alpha '>\alpha $ is finite, for every
$1\le b< \infty $ if $\alpha ' = \infty $, since $E$ is the
$C^{\alpha '}_{\beta }$ fiber bundle, $T^bE$ is the fiber bundle
with the base space $N$. Hence the multiplication $(<f>_{K,\alpha
,\beta }, <g>_{K,\alpha ,\beta }>)\mapsto <f>_{K,\alpha ,\beta
}\circ <g>_{K,\alpha ,\beta }= <f\vee g>_{K,\alpha ,\beta }$ is
continuous in $(S^ME)_{\alpha ,\beta }$ and is of class $C^l_{\beta
}$ with $l=\alpha '-\alpha $ for finite $\alpha '$ and $l=\infty $
for $\alpha ' = \infty$.

\par The topological spaces $E$ and $N$ are of cardinalities
not less than $\sf c$, hence $card [(S^ME)_{\alpha ,\beta }]\ge \sf
c$.

\par {\bf 4. Definition.} The object $(S^ME)_{\alpha ,\beta }$ from Theorem
3 we call the wrap monoid.

\par {\bf 5. Corollary.} {\it Let $\phi : {\bar M}_1\to {\bar M}_2$ be a surjective
$C^{\alpha }_{\beta }$-mapping of $C^{\alpha }_{\beta }$
differentiable spaces over the same Cayley-Dickson algebra ${\cal
A}_r$, $1\le r\le 3$, or a field ${\bf K}={\cal A}_0$, such that
$\phi (s_{1,0,q})=s_{2,0,a(q)}$ for each $q=1,...,k_1$, where $\{
s_{j,0,q}: q=1,...,k_j \} $ are marked points in ${\bar M}_j$,
$j=1,2$, $1\le a\le k_2$, $l_1\le k_2$, $l_1 := card~ \phi (\{
s_{1,0,q}: q=1,...,k_1 \} )$. Then there exists an induced
homomorphism of topological monoids $\phi ^*: (S^{M_2}E)_{\alpha
,\beta } \to (S^{M_1}E)_{\alpha ,\beta }$. If $l_1=k_2$, then $\phi
^*$ is the embedding.}
\par {\bf Proof.} Take $\Xi _1: {\hat M_1}\to {\bar M}_1$ with marked points
$\{ {\hat s}_{1,0,q}: q=1,...,2k_1 \} $ as in \S 1, then take ${\hat
M}_2$ the same ${\hat M}_1$ with additional $2(k_2-l_1)$ marked
points $\{ {\hat s}_{2,0,q} : q=1,...,2k_3 \} $ such that ${\hat
s}_{1,0,q}={\hat s}_{2,0,q}$ for each $q=1,..,k_1$,
$k_3=k_1+k_2-l_1$, then $\phi \circ \Xi _1 := \Xi _2: {\hat M}_2\to
{\bar M}_2$ is the desired mapping inducing the parallel transport
structure from that of $M_1$. Therefore, each ${\hat \gamma }_2:
{\hat M}_2\to N$ induces ${\hat \gamma }_1: {\hat M_1}\to N$ and to
${\bf P}_{{\hat \gamma }_2,u^2}$ there corresponds ${\bf P}_{{\hat
\gamma }_1,u^1}$ with additional conditions in extra marked points,
where $u^1\subset u^2$. The equivalence class $<(E_2,{\bf P}_{{\hat
\gamma }_2,u^2})>_{K,\alpha ,\beta }\in (S^{M_2}E)_{K,\alpha ,\beta
}$ gives the corresponding elements $<(E_1,{\bf P}_{{\hat \gamma
}_1,u^1})>_{K,\alpha ,\beta }\in (S^{M_1}E)_{K,\alpha ,\beta }$,
since $Di^{\alpha }_{\beta ,0}({\hat M}_1, \{ {\hat s}_{0,q}:
q=1,...,2k_2 \} )\subset Di^{\alpha }_{\beta ,0}({\hat M}_1, \{
{\hat s}_{0,q}: q=1,...,2k_3 \} )$. Then $\phi ^*(<(E_2, {\bf
P}_{{\hat \gamma }_2,u^2})\vee (E_1, {\bf P}_{{\hat \eta
}_2,v^2})>_{K,\alpha ,\beta }) =\phi ^*(<(E_2, {\bf P}_{{\hat \gamma
}_2,u^2})>_{K,\alpha ,\beta }) \phi ^*(<(E_1,{\bf P}_{{\hat \eta
}_2,v^2})>_{K,\alpha ,\beta })$, since $f_2\circ \phi (x)$ for each
$x\in {\Xi _1({\hat M}_1\setminus {\hat M}_f})$ coincides with
$f_1(x)$, where $f_j$ corresponds to ${\bf P}_{\gamma _j,y_0\times
e}$ (see also the beginning of \S 3).
\par If $l_1=k_2$, then ${\hat M}_1={\hat M}_2$ and
the family $Di^{\alpha }_{\beta ,0}({\hat M}_1)$ is the same for two
cases, hence $\phi ^*$ is bijective and inevitably $\phi ^*$ is the
embedding.

\par {\bf 6. Theorems.}  {\it {\bf 1.} An alternative
topological Hausdorff group $(W^ME)_{\alpha ,\beta }$ exists
containing the monoid $(S^ME)_{\alpha ,\beta }$ and its group
operation of $C^l_{\beta }$ class is with $l=\alpha ' - \alpha $
($l=\infty $ for $\alpha ' =\infty $). If $M$ and $N$ and $G$ are
separable, then $(W^ME)_{\alpha ,\beta }$ is separable. If $N$ and
$G$ are complete, then $(W^ME)_{\alpha ,\beta }$ is complete. \par
{\bf 2.} If $G$ is associative, then $(W^ME)_{\alpha ,\beta }$ is
associative. If $G$ is commutative, then $(W^ME)_{\alpha ,\beta }$
is commutative. If $G$ is a Lie group, then $(W^ME)_{\alpha ,\beta
}$ is a Lie group. \par {\bf 3.} The group $(W^ME)_{\alpha ,\beta }$
is non-discrete, totally disconnected and non locally compact for
non degenerate $N$. Moreover, if there exist two different sets of
marked points $s_{0,q,j}$ in ${\bar M}_f$, $q=1,...,k$, $j=1,2$,
then two groups $(W^ME)_{\alpha ,\beta ,j}$, defined for $\{
s_{0,q,j}: q=1,...,k \} $ as marked points, are isomorphic.
\par {\bf 4.} The $(W^ME)_{\alpha ,\beta }$ has a structure of an
$C^{\alpha }_{\beta }$-differentiable manifold over ${\cal A}_r$.}

\par {\bf Proof.} If $\gamma \in
C^{\alpha }_{\beta }((M,\{ s_{0,q}: q=1,...,k \} ); (N,y_0))$, then
for $u\in E_{y_0}$ there exists a unique ${\sf h}_q\in G$ such that
${\bf P}_{{\hat \gamma },u}({\hat s}_{0,q+k})=u_q{\sf h}_q$, where
${\sf h}_q=g_q^{-1}g_{q+k}$, $y_0\times g_q = {\bf P}_{{\hat \gamma
},u}({\hat s}_{0,q})$, $g_q\in G$. Due to the equivariance of the
parallel transport structure $h$ depends on $\gamma $ only and we
denote it by ${\sf h}^{(E,{\bf P})}(\gamma )={\sf h}(\gamma )= {\sf
h}$, ${\sf h} = ({\sf h}_1,...,{\sf h}_k)$. The element ${\sf
h}(\gamma )$ is called the holonomy of $\bf P$ along $\gamma $ and
${\sf h}^{(E,{\bf P})}(\gamma )$ depends only on the isomorphism
class of $(E,{\bf P})$ due to the use of the family $Di^{\alpha
}_{\beta ,0}({\hat M})$ and boundary conditions on $\hat \gamma $ at
${\hat s}_{0,q}$ for $q=1,...,2k$.
\par Therefore, ${\sf h}^{(E_1,{\bf P}_1)(E_2,{\bf P}_2)}(\gamma
)={\sf h}^{(E_1,{\bf P}_1)}(\gamma ) {\sf h}^{(E_2,{\bf
P}_2)}(\gamma )\in G^k$, where $G^k$ denotes the direct product of
$k$ copies of the group $G$. Hence for each such $\gamma $ there
exists the homomorphism ${\sf h}(\gamma ): (S^ME)_{\alpha , \beta
}\to G^k$, which induces the homomorphism ${\sf h}: (S^ME)_{\alpha ,
\beta }\to C^0(C^{\alpha }_{\beta }((M,\{ s_{0,q}: q=1,...,k \} );
(N,y_0) ),G^k)$, where $C^0(A,G^k)$ is the space of continuous maps
from a topological space $A$ into $G^k$ and the group structure
$({\sf h} {\sf b})(\gamma )={\sf h}(\gamma ){\sf b}(\gamma )$ (see
also \cite{gajer} for $S^n$).
\par We construct now $(W^MN)_{\alpha ,\beta }$
from $(S^MN)_{\alpha ,\beta }$. In view of Theorem 3 we have the
commutative monoid $(S^MN)_{\alpha ,\beta }$ with the unit and the
cancelation property. Algebraically a group associated with this
monoid is the quotient group $F/\sf B$, where $F$ is the free
commutative group generated by $(S^MN)_{\alpha ,\beta }$, while $\sf
B$ is the minimal closed subgroup in $F$ generated by all elements
of the form $[f+g]-[f]-[g]$ with $f$ and $g\in (S^MN)_{\alpha ,\beta
}$, $[f]$ denotes the element in $F$ corresponding to $f$ (see also
about such abstract Grothendieck's construction in
\cite{langal,swan}).
\par In accordance with Theorem 3 the monoid $(S^MN)_{\alpha ,\beta }$
is the topological $T_1$-space. In view of Theorem 2.3.11 \cite{eng}
the product of $T_1$-spaces is the $T_1$-space. On the other hand,
for a topological group $G$ from the separation axiom $T_0$ it
follows, that $G$ is the Tychonoff space (see Theorems 4.2 and 8.4
in \cite{hew} and also \cite{eng}). The latter means that for a
topological group being $T_0$ or $T_1$ or Hausdorff or Tychonoff is
equivalent. \par At the same time the natural mapping $\eta :
(S^MN)_{\alpha ,\beta }\to (W^MN)_{\alpha ,\beta }$ is injective. We
supply $F$ with the topology inherited from the topology of the
Tychonoff product $(S^MN)_{\alpha ,\beta  }^{\bf Z}$, where each
element $z$ in $F$ has the form $z=\sum_fn_{f,z}[f]$, $n_{f,z}\in
\bf Z$ for each $f\in (S^MN)_{\alpha ,\beta }$,
$\sum_f|n_{f,z}|<\infty $. By the construction $F$ and $F/\sf B$ are
$T_1$-spaces, consequently, $F/\sf B$ is the Tychonoff space. In
particular, $[nf]-n[f]\in \sf B$. We deduce that $\eta $ is the
topological embedding, since $\eta (f+g)=\eta (f)+ \eta (g)$ for
each $f, g \in (S^MN)_{\alpha ,\beta }$, $\eta (e)=e$, since
$(z+B)\in \eta (S^MN)_{\alpha ,\beta }$, when $n_{f,z}\ge 0$ for
each $f$, and inevitably in the general case $z=z^{+}-z^{-}$, where
$(z^{+}+B)$ and $(z^{-}+B)\in \eta (S^MN)_{\alpha ,\beta }$. The
uniform space $(W^ME)_{\alpha , \beta }$ has embedding as the closed
subset into $[(S^ME)_{\alpha ,\beta }]^2$. Thus if $N$ and $G$ are
complete, then $(W^ME)_{\alpha ,\beta }$ is the complete topological
group, since the product of complete uniform spaces is complete (see
Theorem 8.3.9 \cite{eng}) and $(S^ME)_{\alpha ,\beta }$ is complete
by Theorem 3 above.
\par Using plots and $C^{\alpha '}_{\beta }$ transition mappings of
(generalized) charts of $N$ and $E(N,G,\pi ,\Psi )$ and equivalence
classes relative to $Di^{\alpha }_{\beta ,0}(M)$ we get, that
$(W^ME)_{\alpha ,\beta }$ has the structure of the $C^{\alpha
}_{\beta }$-differentiable manifold, since $\alpha ' \ge \alpha $.
\par The rest of the proof and the statements of Theorems 6(1-4)
follows from this and Theorems 3(1-3) and \cite{luannmbp1,lujmslg}.
\par The monoid $(S^ME_U)_{\alpha ,\beta }$ is infinite dimensional over
${\cal A}_r$ due to Theorem 3.3, consequently, $(W^ME_U)_{\alpha
,\beta }$ is infinite dimensional over ${\cal A}_r$, when
$dim_{{\cal A}_r} (X_U) \ge 1$ that is the case, since $N$ is non
degenerate. Thus $(W^ME)_{\alpha ,\beta }$ is non locally compact.

\par {\bf 7. Definition.} The object $(W^ME)_{\alpha ,\beta }
=(W^{M, \{ s_{0,q}: q=1,...,k \} }E; N,G,{\bf P})_{\alpha ,\beta }$
from Theorem 6.1 we call the wrap group.

\par {\bf 8. Corollary.} {\it There exists the group homomorphism
${\sf h}: (W^ME)_{\alpha ,\beta }\to C^0(C^{\alpha }_{\beta }(M,\{
s_{0,q}: q=1,...,k \}; N,y_0),G^k)$.}
\par The {\bf proof} follows from \S 6 and putting ${\sf h}^{f^{-1}}(\gamma )=
({\sf h}^f(\gamma ))^{-1}$.

\par {\bf 9. Corollary.} {\it If $M_1$ and $M_2$ and $\phi $
satisfy conditions of Corollary 5, then there exists a homomorphism
$\phi ^*: (W^{M_2}E)_{\alpha ,\beta }\to (W^{M_1}E)_{\alpha ,\beta
}$. If $l_1=k_2$, then $\phi ^*$ is the embedding.}

\par {\bf 10. Remark.} Each $C^{\alpha '}_{\beta }$ manifold is a
$C^{\alpha '}_{\beta }$ differentiable space. Above differentiable
spaces or manifolds modeled on topological ${\cal A}_r$ vector
spaces were considered. As a particular case of a topological vector
spaces $X$ may be a locally ${\cal A}_r$ convex vector space. It is
well-known that in this case its topology is equivalently
characterized by a family of continuous ultra-pseudo-norms. \par We
recall that a pseudo-norm $v$ on $X$ is called an ultra-pseudo-norm,
if instead of the triangle inequality it satisfies the stronger
condition: $v(x+y)\le \max (v(x), v(y))$ for all $x, y \in X$. A
locally ${\cal A}_r$ convex space $X$ is complete if and only if it
is a projective limit of Banach spaces over ${\cal A}_r$, since $X=
\mbox{}_0Xu_0\oplus ... \oplus \mbox{}_{2^r-1}X u_{2^r-1}$ for each
$1\le r$, where $\mbox{}_0X,...,\mbox{}_{2^r-1}X$ are pairwise
isomorphic locally $\bf K$ convex spaces (see \cite{nari,roo}). For
$r=0$ these spaces are usual $\bf K$-linear spaces. We mention also
that in an ultra-normed space $X$ each two balls $B(X,x,b):= \{ z\in
X: \| z-x \| \le b \} $ and $B(X,y,c)$ either do not intersect or
one of them is contained in another, where $0<b, c <\infty $.
\par Above different classes $C^{\alpha }_{\beta }$ of smoothness were
considered. In particular for $\alpha =0$ this simply reduces to the
class $C^0_{\beta }$ of continuous mappings. For $G= \{ e \} $ there
may $\alpha ' = \alpha $ also be.

\par {\bf 11. Theorem.} {\it For a wrap group $W= (W^ME)_{\alpha ,\beta }$
(see Definition 7 above) there exists a skew product ${\hat W} = W
{\tilde \otimes } W$ which is an $C^l_{\beta }$ alternative group
and there exists a group embedding of $W$ into ${\hat W}$, where
$l=\alpha ' - \alpha $ ($l=\infty $ for $\alpha ' =\infty $),
$E=E(N,G,\pi ,\Psi )$ is a principal $G$-bundle of class $C^{\alpha
'}_{\beta }$ with $\alpha '\ge \alpha \ge 0$. If $G$ is associative,
then $\hat W$ is associative. }

\par {\bf Proof.} Suppose that ${\tilde W}$  is a set of all elements
$(g_1a_1\otimes g_2a_2)\in (W\otimes B)^2$, where $B$ is a free
non-commutative associative group with two generators $a, b$, $ab\ne
ba$, $g_1, g_2\in W$. Take in $\tilde W$ the equivalence relation:
$g_1g_2a\otimes g_2b {\tilde =} ~ g_1e_B\otimes ee_B$, for each
$g_1, g_2\in W$, where $e$ and $e_B$ denote the unit elements in $W$
and in $B$. Define in $\tilde W$ a multiplication by the formula:
\par $(g_1a_1\otimes g_2a_2){\tilde \otimes }
(g_3a_3\otimes g_4a_4) := ((g_1g_3)(a_1a_3)
\otimes (g_4g_2)((a_1^{-1}a_4a_1) a_2)$ \\
for each $g_1, g_2, g_3, g_4\in W$ and every $a_1, a_2, a_3, a_4\in
B$, hence \par $(e\otimes g_1a_1){\tilde \otimes } (e\otimes g_2a_2)
= e\otimes (g_2g_1)(a_2a_1)$, \par $(g_1a_1\otimes e){\tilde \otimes
} (g_2a_2\otimes e) = (g_1g_2)(a_1a_2)\otimes e$,\par
$(g_1a_1\otimes e){\tilde \otimes } (e\otimes g_4a_4) =
g_1a_1\otimes g_4(a_1^{-1}a_4a_1)$, \par $(e\otimes g_4a_4){\tilde
\otimes } (g_1a_1\otimes e) := g_1a_1\otimes g_4a_4$. \\ Thus this
semidirect product $\tilde W$ of groups $(W\otimes B)\otimes ^s
(W\otimes B)$ is non-commutative, since $b^{-1}aba^{-1}\ne e$, where
$e := e\times e_B$, $\otimes ^s$ denotes the semidirect product,
$\otimes $ denotes the direct product.
\par We consider the minimal closed subgroup $A$ in the semidirect product
$\tilde W$ generated by elements $(g_1g_2a\otimes g_2b) {\tilde
\otimes } (g_1e_B\otimes ee_B)^{-1}$, where $B$ is supplied with the
discrete topology and $\tilde W$ is supplied with the product
uniformity. Then put ${\hat W} := {\tilde W}/A =: W{\tilde \otimes
}W$ and denote the multiplication in $\hat W$ as in $\tilde W$. We
get for $W$ the group embedding $\theta : g\mapsto (ge_B\otimes e)$
into ${\hat W}$ and the multiplication $m [(g_1e_B\otimes e),
(g_2e_B\otimes e)] = (g_1e_B\otimes e){\tilde \otimes }
(g_2e_B\otimes e)$.

On the other hand, $(ga_1\otimes e) {\tilde \otimes } (e\otimes
ga_1a_2a_1^{-1}) = ga_1\otimes ga_2 = (e\otimes e) =: {\tilde e}$,
${\hat e} = {\tilde e}A=A$ is the unit element in $\hat W$ and
$(e\otimes ga_1a_2a_1^{-1}) = (ga_1\otimes e)^{-1}$ is the inverse
element of $(ga_1\otimes e)$, where $a_2\in B$ is such that
$(a_1\otimes a_2){\tilde \otimes }A = (e\otimes e){\tilde \otimes
}A=A$ in $\hat W$, $a_1=ea_1$, that is $a_1\otimes a_2 {\tilde =}
e\otimes e$ in $\tilde W$. The preceding formulas mean that $\hat W$
is noncommutative and alternative. \par Moreover, $\hat W$ is the
quotient of a $C^{\alpha }_{\beta }$ differentiable space or a
manifold $W^2$ by the $C^{\alpha }_{\beta }$ equivalence relation
$K_{\alpha ,\beta }$, hence $\hat W$ is the $C^{\alpha }_{\beta }$
differentiable space, since Conditions $(D1-D8)$ of \S 2.4 are
satisfied. The group operation and the inversion in $\hat W$ combine
the product in $W$ and the inversion with the tensor product and the
equivalence relation, hence they are $C^l_{\beta }$ differentiable
with $l=\alpha ' - \alpha $, $l=\infty $ for $\alpha ' = \infty $,
(see \S \S 1.11, 1.12, 1.15 in \cite{souriau} and Theorem 6 above).
\par Then $((g_1\otimes g_2){\tilde \otimes } (g_3\otimes g_4))
{\tilde \otimes }(g_5\otimes g_6) := ((g_1g_3)g_5\otimes
g_6(g_4g_2))$ and \par $(g_1\otimes g_2){\tilde \otimes }
((g_3\otimes g_4)){\tilde
\otimes } (g_5\otimes g_6)) := (g_1(g_3g_5)\otimes (g_6g_4)g_2)$. \\
Therefore, $\hat W$ is alternative, since $W$ is alternative (see
Theorem 6) and $B$ is associative. If $G$ is associative, then $W$
is associative and $\hat W$ is associative.

\par Let us consider the commutator \par $[(g_1a_1\otimes g_2a_2){\tilde {\otimes }}
(g_3a_3\otimes g_4a_4)]{\tilde {\otimes }}[(g_1a_1\otimes
g_2a_2)^{-1} {\tilde {\otimes }}$ \\ $(g_3a_3\otimes g_4a_4)^{-1}] =
\{ ((g_1g_3)(a_1a_3)\otimes (g_4g_2)((a_1^{-1}a_4a_1)a_2)){\tilde
{\otimes }}$ \\ $ [(g_1^{-1}a_1^{-1}\otimes
g_2^{-1}(a_1a_2^{-1}a_1^{-1})) {\tilde {\otimes
}}(g_3^{-1}a_3^{-1}\otimes g_4^{-1}(a_3a_4^{-1}a_3^{-1}))]$
\\ $=((g_1g_3)(a_1a_3)\otimes (g_4g_2)((a_1^{-1}a_4a_1)a_2)
{\tilde {\otimes }} ((g_1^{-1}g_3^{-1})(a_1^{-1}a_3^{-1})\otimes
(g_4^{-1}g_2^{-1})$
\\ $(a_1(a_3a_4^{-1}a_3^{-1})a_1^{-1})(a_1a_2^{-1}a_1^{-1})))
=(((g_1g_3)(g_1^{-1}g_3^{-1}))(a_1a_3a_1^{-1}a_3^{-1})\otimes $ \\
$((g_4^{-1}g_2^{-1})(g_4g_2)) ((a_1a_3)^{-1}
[((a_1a_3)a_4^{-1}(a_1a_3)^{-1})(a_1a_2^{-1}a_1^{-1})](a_1a_3))
((a_1^{-1}a_4a_1)a_2)$.

\par By the definition a minimal closed subgroup generated by products of such elements
is the (topological) commutant ${\tilde W}_c$ of $\tilde W$. The
group $(W^MN)_{\alpha ,\beta }$ is commutative (see Theorem 6$(2)$).
We have $B/B_c = \{ e \} $, the quotient group $G/G_c = G_{ab}$ is
the abelianization of $G$, particularly if $G$ is commutative, then
$G_{ab}=G$, where $G_c$ denotes the (topological) commutant subgroup
of $G$. Therefore, we infer that \par $(W^ME;N,G,{\bf P})_{\alpha ,
\beta }/ [(W^ME;N,G,{\bf P})_{\alpha , \beta }]_c =
(W^ME;N,G_{ab},{\bf P})_{\alpha , \beta }$ \\ and inevitably we get
${\tilde W}/{\tilde W}_c = (W^ME;N,G_{ab},{\bf P})_{\alpha , \beta
}$.

\par Using the equivalence relation in ${\tilde W}$
we deduce that  ${\hat W}/{\hat W}_c = (W^ME;N,G_{ab},{\bf
P})_{\alpha ,\beta }$.

\par {\bf 12. Remark.} We consider the group $B^2\otimes B^2/{\cal E}$, where
an equivalence relation $\cal E$ is induced by that of in $B^2$ as
in $\tilde W$: $(a\otimes b) \approx (e\otimes e)$, the group $B$ is
the same as in \S 11 with two generators $a, b$. Then this gives the
equivalences: $[(a\otimes b)\otimes (a\otimes b)] ~ {\cal E} ~ [(e
\otimes e)\otimes (e \otimes e)] ~ {\cal E} ~ [(e\otimes b) \otimes
(a\otimes e)] \otimes [(e\otimes b) \otimes (a\otimes e)] ~ {\cal E}
~ \{ (e\otimes b) \otimes [(a\otimes e) \otimes (e\otimes b)] \}
\otimes (a\otimes e) ~ {\cal E}  ~ (e\otimes a^{-1}ba) \otimes
(a\otimes e) ~ {\cal E} ~ [(e\otimes ab)\otimes (ba\otimes e)]$ in
$B^2\otimes B^2$, since $B^4$ is the associative group. This implies
the commutativity of the iterated skew product wrap group, when $G$
is commutative, that is $({\hat W}^M({\hat W}^ME)_{\alpha ,\beta
})_{\alpha ,\beta }= (W^M(W^ME)_{\alpha ,\beta })_{\alpha ,\beta }$,
$G=G_{ab}$. In particular, $({\hat W}^M({\hat W}^MN)_{\alpha ,\alpha
})_{\alpha ,\beta }= (W^M(W^MN)_{\alpha ,\beta })_{\alpha ,\beta }$,
where $G= \{ e \} $.

\par {\bf 13. Proposition.} {\it If there exists an
$C^{\alpha '}_{\beta }$-diffeomorphism $\eta : N\to N$ such that
$\eta (y_0) = {y_0}'$, where $\alpha \le \alpha '$ then wrap groups
$(W^ME; y_0)_{\alpha ,\beta }$ and $(W^ME; {y_0}')_{\alpha ,\beta }$
defined with marked points $y_0$ and ${y_0}'$ are $C^l_{\beta
}$-isomorphic as $C^l_{\beta }$-differentiable groups, where
$l=\alpha ' - \alpha $ for finite $\alpha '$, $l=\infty $ for
$\alpha ' = \infty $.}
\par {\bf Proof.} Suppose that $f\in C^{\alpha }_{\beta }({\bar M},E)$, then $\eta \circ \pi
\circ f(s_{0,q})= \eta (y_0)={y_0}'$ for each marked point $s_{0,q}$
in $\bar M$, where $\pi : E\to N$ is the projection, $\pi \circ
f=\gamma $, $\gamma $ is a wrap, that is an $C^{\alpha }_{\beta
}$-mapping from $\bar M$ into $N$ with $\gamma (s_{0,q})=y_0$ for
$q=1,....,k$. The differentiable space $N$ is totally disconnected
together with $E$ and $G$ in accordance with conditions imposed in
Section 2. Consider the $C^{\alpha '}_{\beta }$-diffeomorphism $\eta
\times e$ of the principal bundle $E$. Then $\Theta : C^{\alpha
}_{\beta }({\bar M}, {\cal W})\to C^{\alpha }_{\beta }({\bar
M},{\cal W})$ is the induced isomorphism such that $\pi \circ \Theta
(f) := \eta \circ \pi \circ f: {\bar M}\to N$ and $(\eta \times
e)\circ f = \Theta (f)$ for $f\in C^{\alpha }_{\beta }({\bar M},E)$.
The mapping $\Theta $ is $C^l_{\beta }$ differentiable by $f$, hence
it gives the $C^l_{\beta }$ isomorphism of the considered
$C^l_{\beta }$-differentiable wrap groups (see Theorem 6$(1)$).

\par Ludkovsky S.V. Department of Applied Mathematics MIREA,
av. Vernadsky 78, Moscow 119454
\par sludkowski@mail.ru

\end{document}